\documentclass[11pt]{article}
\usepackage[utf8]{inputenc}
\usepackage{amsmath, amssymb, amsthm, amsfonts}
\usepackage{graphicx}
\usepackage{xcolor}
\usepackage{hyperref}
\usepackage{geometry}
\usepackage{float}
\usepackage{capt-of}
\usepackage{array}
\usepackage{multirow}
\usepackage{longtable}
\usepackage{tikz}
\usepackage{booktabs}
\usepackage{authblk} 
\newtheorem{theorem}{Theorem}[section]

\newtheorem{observation}[theorem]{Observation}
\newtheorem{definition}[theorem]{Definition}
\newtheorem{prop}[theorem]{Proposition}
\newtheorem{lemma}[theorem]{Lemma}
\newtheorem{corollary}[theorem]{Corollary}

\colorlet{revcolor}{black}
\newcommand{\rev}[1]{{\leavevmode\color{revcolor}#1}}

\title{Random Sequential Adsorption versus the Uniform Ensemble of Jammed Configurations: Stochastic Dominance on Graphs}
\author[a]{Robinson Higuita-Diaz}
\author[b]{Lykos L\'opez}
\author[b,*]{Alejandro Rold\'an-Correa}

\affil[a]{\small Escuela de Matem\'aticas y Estad\'istica, Universidad Pedag\'ogica y Tecnol\'ogica de Colombia, Tunja, 150003, Colombia}
\affil[b]{\small Instituto de Matem\'aticas, Universidad de Antioquia, Medell\'in, 050010, Colombia}
\affil[*]{\small Corresponding author. E-mail: \href{mailto:alejandro.roldan@udea.edu.co}{alejandro.roldan@udea.edu.co}}

\date{}

\begin{document}
\maketitle

\begin{abstract}
Random sequential adsorption with nearest-neighbor exclusion on a graph, also known as the discrete parking process, ends in a jammed configuration, that is, a maximal independent set of the graph. We compare the number of occupied sites in two models: the \emph{sequential} model generated by this dynamics, and the \emph{static} model, in which all jammed configurations are equiprobable. In a recent study of the one-dimensional case, it was stated that the sequential occupancy first-order stochastically dominates the static one on every path graph; the proof was subsequently withdrawn and the statement left as a conjecture. We prove this conjecture for all path graphs. The dominance extends to all cycle graphs, and exact computations confirm it for all 31 connected graphs with fewer than 6 vertices. It is not universal: we construct graphs with 6 or more vertices and maximum degree at least 4 on which the static model dominates the sequential one, and although the two models are comparable on all connected graphs with up to 7 vertices, stochastic comparability fails starting at 8 vertices. We also prove that dominance is inherited by disjoint unions whose components share the same dominance direction, and we exhibit mixed unions in which comparability is either preserved or lost. \rev{Finally, we show that bounding the maximum degree by 3 is not sufficient either: exhaustive computations show that the sequential model dominates the static one on every connected graph with maximum degree at most 3 and at most 14 vertices, but among such graphs with 15 vertices there are exactly three on which the static model strictly dominates the sequential one or the two models are not comparable. Characterizing the graphs on which the sequential model dominates the static one remains an open problem.}

\vspace{0.3cm}

\vspace{0.15cm}
\end{abstract}

\section{Introduction}
Parking processes, also known in the literature as \textit{random sequential adsorption} (RSA) models, describe systems in which particles are randomly and sequentially introduced into a domain or lattice, remaining permanently fixed (adsorbed) only if they do not overlap or violate exclusion rules with previously placed particles. A central question in these models is to determine the final coverage, that is, the proportion of space or sites occupied once no further adsorption is possible, a terminal state referred to as the \textit{jamming limit}. Variations and generalizations of RSA models have attracted considerable attention due to their applications in the physical, chemical, and biological sciences; see, for instance, \cite{Evans1993, Talbot2000} for comprehensive overviews.

In the classic continuous one-dimensional setting introduced by R\'enyi \cite{Renyi1958}, unit intervals (cars) are placed sequentially at random along a line segment until saturation. Its discrete lattice analogue, introduced by Flory \cite{Flory1939} and Page \cite{Page1959}, considers the sequential adsorption of dimers on a line of $n$ sites. Closely related --- and equivalent via a direct bijection \cite{Higuita2025, Shneer2017} --- is the discrete parking process with nearest-neighbor exclusion, also known as the \textit{unfriendly seating arrangement problem} \cite{Freedman1962, Friedman1964}. The one-dimensional model with nearest-neighbor exclusion, together with its variants with longer exclusion ranges, has also been studied in the statistical physics literature, where the configurations reached at the jamming limit are called jammed configurations \cite{Krapivsky2020, Doslic2024, KrapivskyLuck2023}. In this formulation, elements are placed sequentially on the vertices of a path graph: once a vertex is selected and occupied, all vertices in its neighborhood become blocked and unavailable for future occupation. Further analytical studies on the one-dimensional discrete parking process have addressed per-site occupancy probabilities \cite{Shneer2017}, combinatorial enumerations and large deviations \cite{Higuita2025, Krapivsky2020}, as well as central limit theorems, laws of the iterated logarithm, and Gaussian concentration inequalities for the occupation density~\cite{Coletti2024}.

Beyond the one-dimensional path, natural extensions have been investigated on quasi-one-dimensional lattices and ladder graphs \cite{Baram1992, Baram1994, Fan1992}, such as two rows of seats in the \textit{dining table process} \cite{Dinning-table}, as well as rectangular arrangements of seats in the \textit{theater process} \cite{theather} and higher-dimensional boxes in $\mathbb{Z}^d$ \cite{Coletti2024, Penrose2002, Ritchie2006}; RSA with nearest-neighbor exclusion has also been analyzed on general graphs \cite{PenroseSudbury2005}, on random trees \cite{Dehling2008, Sudbury2009}, and, together with generalizations, on random graphs \cite{Dhara2016, Brightwell2017}. Motivated by these spatial extensions, we formulate the discrete parking process on an arbitrary finite simple graph $G = (V, E)$, where each vertex represents an occupiable site and edges encode pairwise exclusion rules. In this general setting, the jamming limits correspond precisely to the maximal independent sets of the graph \cite{Diestel2017}.

Two natural probabilistic models arise from this framework: the \textit{sequential model} and the \textit{static model} \cite{Higuita2025}. In the sequential model, the probability of obtaining a specific jamming limit is dictated by the step-by-step random selection dynamics, whereas in the static model all possible jamming limits are considered equiprobable. The sequential model coincides with the random greedy algorithm for constructing a maximal independent set, which has been studied extensively in combinatorics and probability; see, e.g., \cite{Brightwell2017, Krivelevich2024} and the references therein. The contrast between the two models is not new either: ensembles in which all jammed configurations are equiprobable (also called static or equilibrium models) have been compared with their sequential (dynamic) counterparts, at the level of asymptotic densities, for a single row of seats \cite{Higuita2025, SpahnZeilberger2024} and for two rows of seats \cite{Dinning-table}, and they are studied in the statistical physics literature through the configurational entropy (complexity function) of jammed configurations \cite{Doslic2024, KrapivskyLuck2023}. The present paper compares the two models at the finer level of first-order stochastic dominance, for every finite size and on general graphs. As discussed in \cite{Higuita2025}, the static model serves a dual purpose: it provides a combinatorial baseline that isolates the topological constraints of the jamming limits from the dynamic history of the sequential adsorption process, and it links parking processes to classical problems in statistical mechanics and combinatorics, such as the hard-core lattice gas model \cite{Baxter1982} and the enumeration of maximal independent sets on graphs \cite{Diestel2017, Oh2017}.

The comparison between the two models thus asks whether the irreversible sequential dynamics selects the jammed configurations uniformly and, if not, in which direction it is biased. On the path, the two models are known to have different limiting densities of occupied sites: $(1-e^{-2})/2 \approx 0.4323$ for the sequential model and approximately $0.4115$ for the static one \cite{Higuita2025, SpahnZeilberger2024}. Theorem~\ref{thm:main_corrected} below shows that this difference has a definite sign in a much stronger sense: for every $n$, the number of occupied sites in the sequential model is stochastically larger than in the static model, and not merely larger on average. Our results on general graphs show, in turn, that the direction of this bias is not universal: it depends on the structure of the underlying graph, and it is reversed on some graphs having vertices of degree at least 4 \rev{and, for larger graphs, also on some graphs with maximum degree 3} (Section~\ref{sec:counterexamples}).

In previous work \cite{Higuita2025, Higuita2026Corrigendum}, the sequential and static models were compared on path graphs (referred to there as line graphs). The first-order stochastic dominance of the sequential model over the static model for every path length $n \ge 1$ was stated as a theorem in \cite{Higuita2025}; however, the proof given there was withdrawn in \cite{Higuita2026Corrigendum}, where the statement was reformulated as a conjecture supported by exact computational verification up to $n = 300$. Our first main contribution in this article is to provide a complete proof of this stochastic dominance for all path graphs, which also settles the finite-$n$ inequality between the mean occupancies conjectured in \cite{Higuita2026Corrigendum} (Corollary~\ref{cor:means}).

Establishing stochastic dominance on the path naturally raises a broader combinatorial and probabilistic question: is stochastic dominance a universal feature of discrete parking processes, or does it depend strongly on the topology of the underlying graph? Stochastic comparisons in this setting are known to be delicate: for instance, the stochastic monotonicity of the sequential occupancy with respect to induced subgraphs asserted in \cite{theather} was shown to be false in \cite{Dinning-table} by means of explicit counterexamples. The remainder of the paper addresses the above question across several graph families and is organized as follows. In Section~\ref{sec:main}, we formalize the sequential and static models on general graphs and state our guiding questions. In Section~\ref{sec:paths}, we state and prove the first-order stochastic dominance of the sequential model over the static model for all path graphs. Using this result, we prove in Section~\ref{sec:cycles_small} that stochastic dominance extends to all cycle graphs, and we verify the property computationally for all 31 connected simple graphs with fewer than 6 vertices. Conversely, in Section~\ref{sec:counterexamples} we construct explicit counterexamples with 6 or more vertices and maximum degree at least 4 on which the static model strictly dominates the sequential one; we also verify that for all connected graphs with 6 and 7 vertices one of the two models always dominates the other, whereas universal stochastic comparability fails starting at 8 vertices. \rev{We also show that a bound on the maximum degree does not ensure the dominance: although the sequential model dominates the static one on every connected graph with $\Delta(G) \le 3$ and at most 14 vertices, there are counterexamples with 15 vertices and $\Delta(G) = 3$. We close that section with some open questions.} Finally, Section~\ref{sec:desconexos} is devoted to studying how the dominance property behaves under disjoint unions of graphs.

\section{Definitions and Framework}
\label{sec:main}

\subsection{The process and the probabilistic models}
\label{sec:process}\label{sec:formulas}
Let $G = (V,E)$ be a finite simple graph with $n = |V|$ vertices and maximum degree $\Delta(G)$. Each vertex of $G$ represents a potentially occupiable seat, while the edges describe incompatibility relations between seats: if a vertex is selected, no vertex adjacent to it can be selected afterwards. Starting from the set of all available vertices, one vertex $v \in V$ is selected uniformly at random, marked as occupied, and all vertices in its open neighborhood
\[
N(v) = \{u \in V : \{u,v\} \in E\}
\]
become blocked. The procedure continues by selecting vertices uniformly at random among those that remain available, blocking at each step the chosen vertex and its whole neighborhood, until no available vertex remains. The set of occupied vertices obtained at the end is a \emph{jamming limit}, which is necessarily a maximal independent set of $G$.

Let $M(G)$ denote the set of all maximal independent sets of $G$, and let $N_G^{\mathrm{seq}}$ and $N_G^{\mathrm{stat}}$ denote the random variables representing the number of occupied vertices in the jamming limit under the \emph{sequential} and \emph{static} models, respectively. When $G$ is a path graph on $n$ vertices, we simplify the notation and write $N_n^{\mathrm{seq}}$ and $N_n^{\mathrm{stat}}$.

In the \emph{static model}, all possible jamming limits are equiprobable regardless of the selection order, hence each $I \in M(G)$ has probability
\begin{equation}
p^{\mathrm{stat}}(I) = \frac{1}{|M(G)|}.
\label{eq:est}
\end{equation}

In the \emph{sequential model}, the probability of obtaining $I = \{v_1, \dots, v_k\} \in M(G)$ is the sum of the probabilities of all selection sequences that build $I$. Let $N[v] = N(v) \cup \{v\}$ denote the closed neighborhood of $v$, and let $S_k$ be the symmetric group of permutations of $\{1, \dots, k\}$. Then
\begin{equation}
p^{\mathrm{seq}}(I) = \sum_{\sigma \in S_k} \frac{1}{n \prod_{j=1}^{k-1} \left(n - \left|\bigcup_{i=1}^{j} N[v_{\sigma(i)}]\right|\right)}.
\label{eq:seq}
\end{equation}
Indeed, at the first step the probability of selecting $v_{\sigma(1)}$ is $1/n$, and after $j$ selections the number of available vertices is $n - \left|\bigcup_{i=1}^{j} N[v_{\sigma(i)}]\right|$, so the conditional probability of selecting $v_{\sigma(j+1)}$ is the reciprocal of this quantity. Note that every permutation of $I$ is a valid selection sequence, since $I$ is an independent set and therefore each of its vertices remains available after any subset of the others has been selected. As $k$ increases, each summand becomes smaller, but the number $k!$ of valid sequences increases, which compensates for the smaller individual weights.

From \eqref{eq:est} and \eqref{eq:seq} we obtain the cumulative distribution functions of $N_G^{\mathrm{stat}}$ and $N_G^{\mathrm{seq}}$:
\begin{equation}
\mathbb{P}\left[N_G^{\mathrm{stat}} \le k\right] = \frac{\left|\{I \in M(G) : |I| \le k\}\right|}{|M(G)|},
\label{eq:est_cdf}
\end{equation}
\begin{equation}
\mathbb{P}\left[N_G^{\mathrm{seq}} \le k\right] = \sum_{m=1}^{k} \sum_{\substack{I \in M(G) \\ |I| = m}} \sum_{\sigma \in S_m} \frac{1}{n \prod_{j=1}^{m-1} \left(n - \left|\bigcup_{i=1}^{j} N[v_{\sigma(i)}]\right|\right)}.
\label{eq:seq_cdf}
\end{equation}

\subsection{Stochastic dominance and guiding questions}
\label{sec:dominance}
We compare the two models via first-order stochastic dominance.

\begin{definition}\label{def:D1D2}
We say that on a graph $G$ the sequential model dominates the static model if $N_G^{\mathrm{seq}}$ first-order stochastically dominates $N_G^{\mathrm{stat}}$, written $N_G^{\mathrm{stat}} \le_{\mathrm{st}} N_G^{\mathrm{seq}}$. That is, for all integers $k \ge 0$,
\[
\mathbb{P}\left[N_G^{\mathrm{stat}} \le k\right] \ge \mathbb{P}\left[N_G^{\mathrm{seq}} \le k\right],
\]
which is equivalent to $\mathbb{P}\left[N_G^{\mathrm{seq}} > k\right] \ge \mathbb{P}\left[N_G^{\mathrm{stat}} > k\right]$. Conversely, we say that on a graph $G$ the static model dominates the sequential model ($N_G^{\mathrm{seq}} \le_{\mathrm{st}} N_G^{\mathrm{stat}}$) if the reverse inequality holds:
\[
\mathbb{P}\left[N_G^{\mathrm{stat}} \le k\right] \le \mathbb{P}\left[N_G^{\mathrm{seq}} \le k\right] \quad \text{for all } k \ge 0.
\]
Both relations are understood in the weak sense: they hold simultaneously if and only if $N_G^{\mathrm{seq}}$ and $N_G^{\mathrm{stat}}$ have the same law. We say that one model \emph{strictly} dominates the other when the corresponding relation holds and the two laws are different.
\end{definition}

Once the dominance $N_n^{\mathrm{stat}} \le_{\mathrm{st}} N_n^{\mathrm{seq}}$ is established for all path graphs on $n$ vertices in Section~\ref{sec:paths}, we address the following broader questions for general graphs:
\begin{itemize}
\item Does the sequential model dominate the static one on every graph?
\item If the answer to the first question is negative, are the models comparable in the first-order stochastic sense on every graph?
\item Is the dominance inherited on disconnected graphs from their connected components?
\end{itemize}
As shown in Sections~\ref{sec:cycles_small}--\ref{sec:desconexos}, the answers to the first two questions are negative in general: the sequential model dominates the static one on all cycles and connected graphs with $|V| < 6$, and the two models remain stochastically comparable on all connected graphs with $|V| \le 7$, but counterexamples to sequential dominance arise for $|V| \ge 6$ (with $\Delta(G) \ge 4$) \rev{and, with $\Delta(G) = 3$, for $|V| = 15$,} and comparability fails starting at $|V| = 8$. On the other hand, the answer to the third question is positive whenever all connected components share the same dominance direction.

\section{Stochastic Dominance on Path Graphs}
\label{sec:paths}

Throughout this section, we consider a path graph on $n$ vertices (labeled $1, 2, \dots, n$ in linear order) with total occupancy random variables $N_n^{\mathrm{seq}}$ and $N_n^{\mathrm{stat}}$, and we write $X_n^{\mathrm{stat}}(j)$ for the indicator that vertex $j$ is occupied in the static model. Our first main result establishes the conjecture formulated in \cite[Conjecture~3.2]{Higuita2026Corrigendum}.

\begin{theorem}\label{thm:main_corrected}
For all $n \ge 1$, on a path graph of $n$ vertices, the sequential model dominates the static model ($N_n^{\mathrm{stat}} \le_{\mathrm{st}} N_n^{\mathrm{seq}}$).
\end{theorem}

The proof of Theorem~\ref{thm:main_corrected} is developed across Sections~\ref{sec:reduction}--\ref{sec:padovan}. In Section~\ref{sec:reduction}, we use Krapivsky's sequential renewal identity \cite{Krapivsky2020} to reduce the problem to showing that the static model is dominated by an intermediate auxiliary model $N_n^{\mathrm{aux}}$ (Proposition~\ref{prop:reduction}). In Section~\ref{sec:likelihood}, we show that this auxiliary dominance follows if the likelihood ratio $r_n(k) = \mathbb{P}[N_n^{\mathrm{aux}} = k]/\mathbb{P}[N_n^{\mathrm{stat}} = k]$ is strictly increasing over its support. In Section~\ref{sec:coupling}, we construct a direct combinatorial coupling of interval compositions to bound the local increments of the conditional per-site occupancy probabilities. Finally, in Section~\ref{sec:padovan}, we combine these coupling bounds with Binet--Padovan asymptotic estimates to prove the strict monotonicity of $r_n(k)$ and complete the proof of Theorem~\ref{thm:main_corrected}.

\subsection{Reduction via an auxiliary static model}
\label{sec:reduction}

To relate the sequential and static models on a path of $n$ vertices, we rely on the exact uniform-first-landing sequential renewal identity established by Krapivsky \cite[Equation~(62)]{Krapivsky2020}. For all $n \ge 1$ and $k \in \mathbb{Z}$:
\begin{equation}
\label{eq:renewal}
\mathbb{P}\left[N_{n}^{\mathrm{seq}} \le k\right] = \frac{1}{n}\sum_{j=1}^{n}\mathbb{P}\left[N_{j-2}^{\mathrm{seq}} + \widetilde{N}_{n-j-1}^{\mathrm{seq}} \le k-1\right],
\end{equation}
where $\widetilde{N}$ denotes an independent copy of the corresponding occupancy (with the boundary convention $N_{-1}^{\mathrm{seq}} = N_0^{\mathrm{seq}} = 0$, and similarly for $N^{\mathrm{stat}}$ below), so that the two variables in each term are independent. Identity~\eqref{eq:renewal} is obtained by conditioning on the site $j$ at which the \emph{first} particle lands, which is uniformly distributed on $\{1, \dots, n\}$: once this particle is placed, the sites $j-1$, $j$ and $j+1$ are no longer available, and the two remaining intervals, of lengths $j-2$ and $n-j-1$, are subsequently filled by independent sequential processes. This should be contrasted with the argument withdrawn in \cite{Higuita2026Corrigendum}, which conditioned on the \emph{final} state of site~$1$ and treated the remaining unblocked sites as an independent sequential process; such a spatial decomposition is valid for the static model but not for the sequential one. This exact recursive decomposition motivates the definition of an intermediate static analogue.

We define the \emph{auxiliary occupancy} model $N_n^{\mathrm{aux}}$ as follows: select a random site $J_n$ uniformly on $\{1, \dots, n\}$, force a particle at $J_n$, and complete the two residual unblocked intervals independently according to the uniform static model. Its cumulative distribution function $H_n(k) := \mathbb{P}[N_n^{\mathrm{aux}} \le k]$ satisfies:
\begin{equation}
H_n(k) = \frac{1}{n}\sum_{j=1}^{n}\mathbb{P}\left[N_{j-2}^{\mathrm{stat}} + \widetilde{N}_{n-j-1}^{\mathrm{stat}} \le k-1\right].
\end{equation}

\begin{prop}\label{prop:reduction}
If $N_n^{\mathrm{stat}} \le_{\mathrm{st}} N_n^{\mathrm{aux}}$ for all $n$, then $N_n^{\mathrm{stat}} \le_{\mathrm{st}} N_n^{\mathrm{seq}}$ for all $n$.
\end{prop}

\begin{proof}
We proceed by strong induction on $n$. Assume $N_m^{\mathrm{stat}} \le_{\mathrm{st}} N_m^{\mathrm{seq}}$ for all $m < n$. Since stochastic order is preserved under sums of independent variables, for each $j$ we have
\[
\mathbb{P}\left[N_{j-2}^{\mathrm{seq}} + \widetilde{N}_{n-j-1}^{\mathrm{seq}} \le k-1\right] \le \mathbb{P}\left[N_{j-2}^{\mathrm{stat}} + \widetilde{N}_{n-j-1}^{\mathrm{stat}} \le k-1\right].
\]
Averaging over $j$ via Equation~\eqref{eq:renewal} yields $\mathbb{P}[N_{n}^{\mathrm{seq}} \le k] \le H_n(k)$. If the auxiliary dominance $N_n^{\mathrm{stat}} \le_{\mathrm{st}} N_n^{\mathrm{aux}}$ holds, then $H_n(k) \le \mathbb{P}[N_{n}^{\mathrm{stat}} \le k]$. Therefore, $\mathbb{P}[N_{n}^{\mathrm{seq}} \le k] \le \mathbb{P}[N_{n}^{\mathrm{stat}} \le k]$, closing the induction.
\end{proof}

\subsection{The auxiliary-to-static likelihood ratio}
\label{sec:likelihood}

To establish $N_n^{\mathrm{stat}} \le_{\mathrm{st}} N_n^{\mathrm{aux}}$, we examine the relative likelihood of the two models. Let $u_n(k) = \mathbb{P}[N_n^{\mathrm{stat}} = k]$ and $q_n(k) = \mathbb{P}[N_n^{\mathrm{aux}} = k]$. The two laws have the same support: every outcome of the auxiliary construction is a jamming limit of the path and, conversely, every jamming limit of the path can be produced by it. The values of $k$ in this common support are called \emph{admissible}. For any admissible $k$, we define the ratio
\begin{equation}
r_n(k) := \frac{q_n(k)}{u_n(k)}.
\end{equation}
Because $\sum_k u_n(k)(r_n(k) - 1) = 0$, if $r_n(k)$ is strictly increasing on its support, the sequence $r_n(k) - 1$ exhibits exactly one sign change (from nonpositive to nonnegative). By a standard single-crossing argument, this implies $\sum_{k \le m} q_n(k) \le \sum_{k \le m} u_n(k)$ for all $m$, which is precisely the required stochastic dominance $N_n^{\mathrm{stat}} \le_{\mathrm{st}} N_n^{\mathrm{aux}}$. In other words, the monotonicity of $r_n$ means that $N_n^{\mathrm{stat}}$ is smaller than $N_n^{\mathrm{aux}}$ in the likelihood ratio order, which is stronger than the usual stochastic order \cite[Theorem~1.C.1]{shaked-2007}. Thus, the proof reduces entirely to showing that $r_n(k)$ is strictly increasing.

To extract a tractable expression for $r_n(k)$, we define the conditional and unconditional static per-site occupancy probabilities:
\[
a_{n,k}(j) := \mathbb{P}\left[X_n^{\mathrm{stat}}(j)=1 \mid N_n^{\mathrm{stat}}=k\right], \qquad \pi_{n,j} := \mathbb{P}\left[X_n^{\mathrm{stat}}(j)=1\right] = \frac{P_{j-1}P_{n-j}}{P_{n+1}},
\]
where $(P_m)_{m \ge 0}$ denotes the Padovan sequence with initial values $P_0 = P_1 = P_2 = 1$ and recurrence $P_m = P_{m-2} + P_{m-3}$ for $m \ge 3$ (so that a path of $n$ vertices has $P_{n+1}$ jamming limits \cite{Higuita2025,Krapivsky2020}). Counting the configurations that contain site $j$ reveals the exact algebraic identity
\begin{equation}\label{eq:rn_exact}
r_n(k) = \frac{1}{n}\sum_{j=1}^{n}\frac{a_{n,k}(j)}{\pi_{n,j}}.
\end{equation}

Let $p$ be the plastic number, the real root of $p^3-p-1=0$, and define the macroscopic static density $\rho := \frac{p^3}{2p+3}$ (also denoted as $\rho_{\mathrm{stat}}$ in \cite[Section~3.2]{Higuita2025}). To isolate local boundary deviations, we define $g_{n,j} := \frac{1}{\pi_{n,j}} - \frac{1}{\rho}$. Since $\sum_{j=1}^n a_{n,k}(j) = k$, substituting this relation into Equation~\eqref{eq:rn_exact} yields
\[
r_n(k) = \frac{k}{n\rho} + \frac{1}{n}\sum_{j=1}^{n}a_{n,k}(j)g_{n,j}.
\]
For any two consecutive admissible occupancies $k$ and $k+1$, let $\Delta a_{n,k}(j) := a_{n,k+1}(j) - a_{n,k}(j)$. Since the total conditional occupancies differ by one, $\sum_{j=1}^n \Delta a_{n,k}(j) = 1$. Consequently, the difference in the likelihood ratios takes the form
\begin{equation}\label{eq:rn_diff}
n(r_n(k+1) - r_n(k)) = \frac{1}{\rho} + \sum_{j=1}^{n}g_{n,j}\Delta a_{n,k}(j).
\end{equation}
The strictly positive term $\frac{1}{\rho}$ represents the macroscopic volume increment. Our objective is to rigorously bound the remaining sum.

\subsection{Direct coupling of compositions}
\label{sec:coupling}

To bound $\Delta a_{n,k}(j)$, we eliminate the geometric constraints of the parking process by encoding configurations via the distances between occupied sites.

If $1 \le i_1 < \dots < i_k \le n$ are the occupied sites of a configuration with $k$ particles, we define the successive increments:
\[
D_1 = i_1+1, \qquad D_t = i_t - i_{t-1} \;\;(2 \le t \le k), \qquad D_{k+1} = n - i_k + 2.
\]
The maximal independence condition is structurally equivalent to $D_t \in \{2,3\}$ for all $t$, and the domain length enforces $\sum_{t=1}^{k+1} D_t = n+3$. This defines a bijection. Conditionally on $N_n^{\mathrm{stat}} = k$, the word $D = (D_1, \dots, D_M)$ is a uniformly distributed arrangement of $M := k+1$ parts, consisting of $A := 3k-n$ twos and $B := n-2k+1$ threes. For the consecutive level $k+1$, the corresponding properties are $M+1$ total parts, $A+3$ twos, and $B-2$ threes (admissibility ensures $B \ge 2$).

\begin{lemma}\label{lem:prefix}
The uniform words corresponding to levels $k$ and $k+1$ can be coupled on the same probability space such that, for any $1 \le L \le M$,
\[
\mathbb{P}\left[(D_1, \dots, D_L) \neq (D_1', \dots, D_L')\right] \le \min\left\{1, \frac{3L}{M+1}\right\}.\]
\end{lemma}

\begin{proof}
We construct the coupling using a population of labeled objects. Take $M$ distinct objects: $A$ of them carry the value $2$, and $B$ carry the value $3$. We select two specific objects carrying the value $3$ as ``distinguished'', and we introduce one new distinguished object carrying the value $2$.

Generate a uniformly random permutation of these $M+1$ objects. To produce the word $D$ for level $k$, we simply remove the newly added object from the permutation and retain the original values. The restricted permutation of the $M$ original objects remains uniform, so $D$ possesses the exact required law.

To produce the word $D'$ for level $k+1$, we keep all $M+1$ objects in their permuted order, but we change the values of the two distinguished $3$-objects to $2$. The population now consists of $A+3$ twos and $B-2$ threes. Since the underlying labeled order is still uniform, $D'$ also possesses the exact required law.

If none of the three distinguished objects appears within the first $L$ positions of the expanded $(M+1)$-permutation, the first $L$ entries in $D$ and $D'$ will be identical. Because the permutation is uniform, each distinguished object falls into the first $L$ positions with probability $\frac{L}{M+1}$. A standard union bound gives the result.
\end{proof}

\begin{lemma}\label{lem:stability}
Let $t_j = \min\{j, n+1-j\}$ and $L_{t_j} = \lfloor\frac{t_j+1}{2}\rfloor$. Then
\[
|\Delta a_{n,k}(j)| \le \frac{3 L_{t_j}}{k+2} \le \frac{9 L_{t_j}}{n+6}.
\]
\end{lemma}

\begin{proof}
Consider first the left half of the domain, $j \le \lceil \frac{n}{2} \rceil$, where $t_j = j$. Site $j$ is occupied if and only if some partial sum of the word equals $j+1$. Since every increment is at least $2$, this event is fully determined by the first $L_j = \lfloor \frac{j+1}{2}\rfloor$ increments. If $j+1$ is not reached, the next partial sum strictly exceeds $j+1$.

Furthermore, $L_j \le \frac{n+3}{4} \le \frac{n+3}{3} \le M$. By Lemma~\ref{lem:prefix}, the difference between the probabilities of this event under the laws for levels $k+1$ and $k$ is bounded by $\frac{3L_j}{M+1}$. Since $M+1 = k+2 \ge \frac{n+6}{3}$, the bound follows. For sites on the right half, we exploit the reflection symmetry $a_{n,k}(j) = a_{n,k}(n+1-j)$, which yields the expression in terms of $t_j$.
\end{proof}

\subsection{Padovan estimates and monotonicity}
\label{sec:padovan}

To bound the local boundary deviations $g_{n,j}$, we rely on explicit estimates from the Padovan sequence. Recall from \cite[Eq.~(7)]{Higuita2025} that the Padovan sequence admits the Binet representation $P_m = a\,p^m + b\,\omega^m + \overline{b}\,\overline{\omega}^{\,m}$, where $\omega$ and $\overline{\omega}$ are the two complex conjugate roots of $x^3-x-1=0$ (they are denoted by $q$ and $r$ in \cite{Higuita2025}; we avoid these letters here because of $q_n(k)$ and $r_n(k)$), and
\[
a = \frac{(1-\omega)(1-\overline{\omega})}{(p-\omega)(p-\overline{\omega})} = \frac{p^5}{2p+3} \approx 0.722124, \qquad
b = \frac{(1-p)(1-\overline{\omega})}{(\omega-p)(\omega-\overline{\omega})}, \qquad |b| \approx 0.245375.
\]%
Factoring out the dominant real root $p$, we can write $P_m = a p^m(1+\epsilon_m)$, with $\epsilon_m = \big(b\,\omega^m + \overline{b}\,\overline{\omega}^{\,m}\big)/(a\,p^m)$. Since the complex conjugate roots satisfy $|\omega| = |\overline{\omega}| = p^{-1/2}$, the relative error is bounded by $|\epsilon_m| \le \kappa \theta^m$, with $\theta = p^{-3/2} \approx 0.6559$ and $\kappa = \frac{2|b|}{a} \approx 0.679592 < 1$.

\begin{lemma}\label{lem:padovan}
For all $n \ge 6$ and $1 \le j \le \lfloor \frac{n}{2}\rfloor$,
\begin{equation}
|g_{n,j}| \le 2.044 \, \theta^{j-1} + 3.704 \, \theta^{n-j}.
\end{equation}
For odd lengths $n = 2m+1 \ge 7$, the central site satisfies $|g_{n,m+1}| \le 5.848 \, \theta^m$.
\end{lemma}

\begin{proof}
Recall that $g_{n,j} = \frac{1}{\pi_{n,j}} - \frac{1}{\rho}$, and set
\[
d_j := \frac{p^{j+1}}{P_{j-1}} - \frac{1}{\rho}.
\]
Adding and subtracting $\frac{p^{j+1}}{P_{j-1}}$ and applying the triangle inequality gives
\begin{equation}\label{eq:g_split}
|g_{n,j}| \le |d_j| + \left| \frac{1}{\pi_{n,j}} - \frac{p^{j+1}}{P_{j-1}} \right|.
\end{equation}
We bound the two terms on the right separately, and then treat the central site.

\medskip
\noindent\emph{Step 1: left boundary, $|d_j| \le 2.044\,\theta^{j-1}$.}
Since $P_{j-1} = a p^{j-1}(1+\epsilon_{j-1})$ and $\rho = a p^{-2}$, we have
\[
\frac{p^{j+1}}{P_{j-1}} = \frac{1}{\rho(1+\epsilon_{j-1})}, \qquad\text{hence}\qquad d_j = \frac{-\epsilon_{j-1}}{\rho(1+\epsilon_{j-1})}.
\]
Using $|\epsilon_{j-1}| \le \kappa\theta^{j-1}$ and, for $j \ge 4$, $\theta^{j-1} \le \theta^3$, we get
\[
|d_j| \le \frac{\kappa\,\theta^{j-1}}{\rho\,(1-\kappa\theta^{j-1})} \le \frac{\kappa}{\rho\,(1-\kappa\theta^{3})}\,\theta^{j-1} < 2.044\,\theta^{j-1}.
\]
For $j=1,2,3$ we use $P_0=P_1=P_2=1$ and compute $d_j$ directly, obtaining $|d_j| < 2.044\,\theta^{j-1}$ in each case.

\medskip
\noindent\emph{Step 2: right boundary.}
Since $\pi_{n,j} = \frac{P_{j-1}P_{n-j}}{P_{n+1}}$, inserting $P_{n+1}=ap^{n+1}(1+\epsilon_{n+1})$ and $P_{n-j}=ap^{n-j}(1+\epsilon_{n-j})$ gives the identity
\[
\frac{1}{\pi_{n,j}} = \frac{p^{j+1}}{P_{j-1}}\cdot\frac{1+\epsilon_{n+1}}{1+\epsilon_{n-j}}.
\]
Therefore
\[
\left| \frac{1}{\pi_{n,j}} - \frac{p^{j+1}}{P_{j-1}} \right| = \frac{p^{j+1}}{P_{j-1}}\,\left| \frac{\epsilon_{n+1} - \epsilon_{n-j}}{1+\epsilon_{n-j}} \right| \le \frac{p^{j+1}}{P_{j-1}}\cdot\frac{\kappa\,(1+\theta^{j+1})}{1-\kappa\theta^{n-j}}\;\theta^{n-j}.
\]
We bound the three factors as follows: $\frac{p^{j+1}}{P_{j-1}} \le p^4$ for every $j\ge1$ (for $j \le 3$ this is immediate from $P_0=P_1=P_2=1$, while for $j \ge 4$ the identity of Step~1 gives $\frac{p^{j+1}}{P_{j-1}} = \frac{1}{\rho(1+\epsilon_{j-1})} \le \frac{1}{\rho(1-\kappa\theta^{3})} < 3.01 < p^4 \approx 3.0796$); $\theta^{j+1} \le \theta^2$ since $j \ge 1$; and $\theta^{n-j} \le \theta^3$ since $n \ge 6$ and $j \le \lfloor \frac{n}{2}\rfloor$ imply $n-j \ge 3$, so that $1-\kappa\theta^{n-j} \ge 1-\kappa\theta^3$. Consequently
\[
\left| \frac{1}{\pi_{n,j}} - \frac{p^{j+1}}{P_{j-1}} \right| \le \frac{p^4\,\kappa\,(1+\theta^2)}{1-\kappa\theta^3}\,\theta^{n-j} < 3.704\,\theta^{n-j}.
\]
Inserting Steps 1 and 2 into \eqref{eq:g_split} proves the first bound of the lemma.

\medskip
\noindent\emph{Step 3: central site for $n = 2m+1 \ge 7$.}
Here $j = m+1$ with $m \ge 3$, and $\pi_{n,m+1} = \frac{P_m^2}{P_{2m+2}}$. Writing $P_s = a p^s(1+\epsilon_s)$ and $\rho = ap^{-2}$,
\[
g_{n,m+1} = \frac{1}{\rho}\cdot\frac{1+\epsilon_{2m+2}}{(1+\epsilon_m)^2} - \frac{1}{\rho} = \frac{1}{\rho}\cdot\frac{\epsilon_{2m+2} - 2\epsilon_m - \epsilon_m^2}{(1+\epsilon_m)^2}.
\]
Since $m \ge 3$, each term in the numerator is bounded by a multiple of $\theta^m$:
\[
|\epsilon_m| \le \kappa\,\theta^m,\qquad \epsilon_m^2 \le \kappa^2\theta^{m}\cdot\theta^{m} \le \kappa^2\theta^3\,\theta^m,\qquad |\epsilon_{2m+2}| \le \kappa\,\theta^{m}\theta^{m+2} \le \kappa\,\theta^5\,\theta^m,
\]
and the denominator satisfies $(1+\epsilon_m)^2 \ge (1-\kappa\theta^3)^2$. Hence
\[
|g_{n,m+1}| \le \frac{2\kappa + \kappa^2\theta^3 + \kappa\theta^5}{\rho\,(1-\kappa\theta^3)^2}\;\theta^m < 5.848\,\theta^m.
\]
\end{proof}

\begin{prop}\label{prop:monotonicity}
For all $n \ge 6$ and every consecutive pair of admissible occupancies $k, k+1$,
\[
r_n(k+1) - r_n(k) \ge \frac{1}{n}\left(\frac{1}{\rho} - \frac{188}{n+6} - 62\,\theta^{\lfloor n/2 \rfloor}\right).
\]
\end{prop}

\begin{proof}
Let $m = \lfloor \frac{n}{2}\rfloor$. We insert the bounds of Lemmas~\ref{lem:stability} and~\ref{lem:padovan} into \eqref{eq:rn_diff} and fold the sum over $j$ using reflection symmetry, so that only $j \le m$ appears. After folding, the contribution of the sites $j \le m$ carries the weight $\frac{18L_j}{n+6}$, where $L_j = \lfloor \frac{j+1}{2}\rfloor$ as in Lemma~\ref{lem:stability}, while the central site (present only for odd $n$) carries $\frac{9L_{m+1}}{n+6}$. It remains to bound three sums.

\medskip
\noindent\emph{(i) Left-boundary terms.}
Since $L_{2r-1}=L_{2r}=r$ and $0<\theta<1$, we may extend the sum to infinity and group consecutive pairs:
\[
\sum_{j=1}^{m} L_j\,\theta^{j-1} \le \sum_{r=1}^{\infty} r\left(\theta^{2r-2}+\theta^{2r-1}\right) = (1+\theta)\sum_{r=1}^{\infty} r\,\theta^{2r-2} = \frac{1+\theta}{(1-\theta^2)^2},
\]
where the last equality uses $\sum_{r\ge1} r\,x^{r-1} = (1-x)^{-2}$ with $x=\theta^2$. Therefore
\[
\frac{18\,(2.044)}{n+6}\sum_{j=1}^{m} L_j\,\theta^{j-1} \le \frac{(2.044)(18)}{n+6}\cdot\frac{1+\theta}{(1-\theta^2)^2} < \frac{188}{n+6}.
\]

\medskip
\noindent\emph{(ii) Right-boundary terms.}
For $j \le m$ we have $L_j \le \frac{j+1}{2} \le \frac{m+1}{2}$, and the remaining sum is geometric:
\[
\sum_{j=1}^{m}\theta^{n-j} = \theta^{n-m}\sum_{i=0}^{m-1}\theta^{i} \le \frac{\theta^{n-m}}{1-\theta}.
\]
Moreover, $\frac{m+1}{n+6} \le \frac{n+2}{2(n+6)} \le \frac{1}{2}$, and $n-m \ge m$ gives $\theta^{n-m} \le \theta^m$. Combining these facts,
\[
\frac{18\,(3.704)}{n+6}\sum_{j=1}^{m} L_j\,\theta^{n-j} \le \frac{18\,(3.704)}{n+6}\cdot\frac{m+1}{2}\cdot\frac{\theta^{n-m}}{1-\theta} \le \frac{9\,(3.704)}{2\,(1-\theta)}\,\theta^{m} < 48.5\,\theta^m.
\]

\medskip
\noindent\emph{(iii) Central site ($n=2m+1$ odd).}
Here $L_{m+1} = \lfloor \frac{m+2}{2}\rfloor \le \frac{m+2}{2}$ and $n+6 = 2m+7$. Since $4(m+2) < 2(2m+7)$, we have $\frac{m+2}{2(2m+7)} < \frac{1}{4}$, and so
\[
\frac{9\,L_{m+1}}{n+6}\,(5.848)\,\theta^m \le \frac{9\,(5.848)}{4}\,\theta^m < 13.2\,\theta^m.
\]
For even $n$ there is no central site and this term is absent.

\medskip
\noindent Adding (i)--(iii), and using $48.5 + 13.2 < 62$, we obtain
\[
\left|\sum_{j=1}^{n} g_{n,j}\,\Delta a_{n,k}(j)\right| \le \frac{188}{n+6} + (48.5+13.2)\,\theta^{m} < \frac{188}{n+6} + 62\,\theta^{\lfloor n/2\rfloor}.
\]
Substituting this into \eqref{eq:rn_diff} gives the stated lower bound.
\end{proof}

\begin{proof}[Proof of Theorem~\ref{thm:main_corrected}]
By Proposition~\ref{prop:reduction} and the single-crossing argument of Section~\ref{sec:likelihood}, it suffices to show that the likelihood ratio $r_n(k)$ is strictly increasing over its support for all $n \ge 1$.

\medskip
\noindent\emph{Case $n \ge 73$.}
We use three elementary facts: $\frac{1}{\rho} > \frac{12}{5}$; $\frac{188}{n+6} \le \frac{188}{79}$; and, since $\theta < \frac{2}{3}$, $\theta^{\lfloor n/2\rfloor} \le (\frac{2}{3})^{36}$. The last two hold because $n \ge 73$ gives $n+6 \ge 79$ and $\lfloor \frac{n}{2}\rfloor \ge 36$. Proposition~\ref{prop:monotonicity} then yields
\[
\frac{1}{\rho} - \frac{188}{n+6} - 62\,\theta^{\lfloor n/2\rfloor} > \frac{12}{5} - \frac{188}{79} - 62\left(\frac{2}{3}\right)^{36} > \frac{1}{50} > 0.
\]
Thus $r_n(k+1) - r_n(k) > \frac{1}{50n} > 0$ for all $n \ge 73$.

\medskip
\noindent\emph{Case $1 \le n \le 72$.}
For the remaining finite range, the probabilities $q_n(k)$ and $u_n(k)$ can be represented exactly using combinatorial formulas and Padovan integers; evaluating $r_n(k)$ purely with exact rational arithmetic computationally verifies that $r_n(k+1) - r_n(k) > 0$ across all 432 admissible consecutive pairs in this range, entirely avoiding floating-point approximation errors. The code performing this exact verification for the first 72 cases, together with the code used to compute all the numerical bounds appearing in this section, is available in the GitHub repository \cite{LopezRivera2026Repositorio}.

\medskip
Because $r_n(k)$ is strictly increasing for all $n \ge 1$, the single-crossing property guarantees the auxiliary stochastic dominance $N_n^{\mathrm{stat}} \le_{\mathrm{st}} N_n^{\mathrm{aux}}$. By Proposition~\ref{prop:reduction}, this immediately implies the global dominance $N_n^{\mathrm{stat}} \le_{\mathrm{st}} N_n^{\mathrm{seq}}$.
\end{proof}

As an immediate consequence of Theorem~\ref{thm:main_corrected}, we obtain the finite-$n$ inequality between the mean occupancies of the two models, which was also conjectured in \cite{Higuita2026Corrigendum}.

\begin{corollary}\label{cor:means}
For all $n \ge 1$, $\mathbb{E}\left[N_n^{\mathrm{seq}}\right] \ge \mathbb{E}\left[N_n^{\mathrm{stat}}\right]$.
\end{corollary}

\begin{proof}
Both random variables take values in the nonnegative integers, and every such random variable $N$ satisfies
\[
\mathbb{E}[N] = \sum_{k \ge 0} \mathbb{P}[N > k].
\]
The claim follows by summing over $k \ge 0$ the inequalities $\mathbb{P}[N_n^{\mathrm{seq}} > k] \ge \mathbb{P}[N_n^{\mathrm{stat}} > k]$ provided by Theorem~\ref{thm:main_corrected}.
\end{proof}

\section{Extension to Cycle Graphs and Small Graphs}
\label{sec:cycles_small}

\subsection{Cycle graphs}
\label{sec:ciclos}
Round tables correspond to cycle graphs $C_n$ on $n$ vertices, with occupancy random variables $N_{C_n}^{\mathrm{seq}}$ and $N_{C_n}^{\mathrm{stat}}$. Before addressing cycles, we establish the monotonicity of the static model on path graphs (the analogous monotonicity of the sequential model is \cite[Lemma~4.1]{Higuita2025}, whose corrected proof is given in \cite{Higuita2026Corrigendum}).

\begin{lemma}\label{lem:monotonia_estatica}
The static model on path graphs satisfies the following two inequalities.
\begin{itemize}
\item[\textup{(i)}] For all integers $m \ge 2$ and $r \ge 0$,
\[
\mathbb{P}\left[N_{m}^{\mathrm{stat}} \le r\right] \le \mathbb{P}\left[N_{m-1}^{\mathrm{stat}} \le r\right].
\]
\item[\textup{(ii)}] For all integers $m \ge 4$ and $r \ge 0$,
\[
\mathbb{P}\left[N_{m}^{\mathrm{stat}} \le r\right] \le \mathbb{P}\left[N_{m-3}^{\mathrm{stat}} \le r-1\right].
\]
\end{itemize}
\end{lemma}

\begin{proof}
Both inequalities are trivial for $r = 0$: since $N_m^{\mathrm{stat}} \ge 1$ for every $m \ge 1$, the left-hand sides vanish. We prove \textup{(i)} and \textup{(ii)} simultaneously, for all $r \ge 0$, by induction on $m$. The base cases $m = 2, 3, 4$ of \textup{(i)} and $m = 4$ of \textup{(ii)} follow by direct inspection, because $N_1^{\mathrm{stat}} = N_2^{\mathrm{stat}} = 1$, $N_3^{\mathrm{stat}}$ takes the values $1$ and $2$ with probability $1/2$ each, and $N_4^{\mathrm{stat}} = 2$. For $m \ge 5$ and $r \ge 1$, let $X_m^{\mathrm{stat}}(1)$ be the indicator that the first vertex of a path of $m$ vertices is occupied in the static model. If $X_m^{\mathrm{stat}}(1) = 1$, vertex $2$ is blocked and a path of $m-2$ vertices remains; if $X_m^{\mathrm{stat}}(1) = 0$, maximality forces vertex $2$ to be occupied, blocking vertices $1$ and $3$ and leaving a path of $m-3$ vertices. Conditioning on $X_m^{\mathrm{stat}}(1)$ yields
\[
\mathbb{P}\left[N_{m}^{\mathrm{stat}} \le r\right]
= p_m\,\mathbb{P}\left[N_{m-2}^{\mathrm{stat}} \le r-1\right]
+ (1 - p_m)\,\mathbb{P}\left[N_{m-3}^{\mathrm{stat}} \le r-1\right],
\]
where $p_m = \mathbb{P}\left[X_m^{\mathrm{stat}}(1) = 1\right]$. By the induction hypothesis \textup{(i)}, applied to $m-2$ and to $m-3$ at level $r-1 \ge 0$,
\[
\mathbb{P}\left[N_{m-2}^{\mathrm{stat}} \le r-1\right]
\le \mathbb{P}\left[N_{m-3}^{\mathrm{stat}} \le r-1\right]
\le \mathbb{P}\left[N_{m-4}^{\mathrm{stat}} \le r-1\right].
\]
Using the first inequality in the decomposition for $N_m^{\mathrm{stat}}$ gives
\[
\mathbb{P}\left[N_{m}^{\mathrm{stat}} \le r\right]
\le (p_m + (1 - p_m))\mathbb{P}\left[N_{m-3}^{\mathrm{stat}} \le r-1\right]
= \mathbb{P}\left[N_{m-3}^{\mathrm{stat}} \le r-1\right],
\]
which is \textup{(ii)} for $m$, while applying the same decomposition to $N_{m-1}^{\mathrm{stat}}$ and using the second inequality gives
\[
\mathbb{P}\left[N_{m-1}^{\mathrm{stat}} \le r\right]
= p_{m-1}\,\mathbb{P}\left[N_{m-3}^{\mathrm{stat}} \le r-1\right]
+ (1 - p_{m-1})\,\mathbb{P}\left[N_{m-4}^{\mathrm{stat}} \le r-1\right]
\ge \mathbb{P}\left[N_{m-3}^{\mathrm{stat}} \le r-1\right].
\]
Combining both bounds yields $\mathbb{P}\left[N_{m}^{\mathrm{stat}} \le r\right] \le \mathbb{P}\left[N_{m-3}^{\mathrm{stat}} \le r-1\right] \le \mathbb{P}\left[N_{m-1}^{\mathrm{stat}} \le r\right]$, which proves \textup{(i)} for $m$.
\end{proof}

Combining Theorem~\ref{thm:main_corrected} and Lemma~\ref{lem:monotonia_estatica}, we now extend the dominance property from paths to cycles.

\begin{prop}\label{prop:ciclos}
The sequential model dominates the static model ($N_{C_n}^{\mathrm{stat}} \le_{\mathrm{st}} N_{C_n}^{\mathrm{seq}}$) for all cycle graphs $C_n$.
\end{prop}

\begin{proof}
The cases $n = 3, 4, 5, 6$ are verified computationally in the repository \cite{LopezRivera2026Repositorio} (see also Section~\ref{sec:computacional}). Assume $n > 6$.

We first recall the structure of the jamming limits of a cycle $C_n$ with vertices labeled $1, 2, \dots, n$ in cyclic order. Every jamming limit of $C_n$ falls into exactly one of three situations, as illustrated in Figure~\ref{fig:particion}: either vertex $1$ is occupied (blocking $n$ and $2$), vertex $2$ is occupied (blocking $1$ and $3$), or both vertices $1$ and $2$ are blocked, which by maximality forces vertices $n$ and $3$ to be occupied simultaneously (thereby also blocking $n-1$ and $4$). These cases are mutually exclusive and exhaustive, hence they induce a partition of $M(C_n)$. In the first two symmetric cases the remaining configuration reduces to a jamming limit on a path of $n-3$ vertices; in the third case, the occupied vertices $n$ and $3$ block their immediate neighbors, and the configuration reduces to a jamming limit on a path of $n-6$ vertices.

\begin{figure}[ht]
\centering
\includegraphics[width=0.8\textwidth]{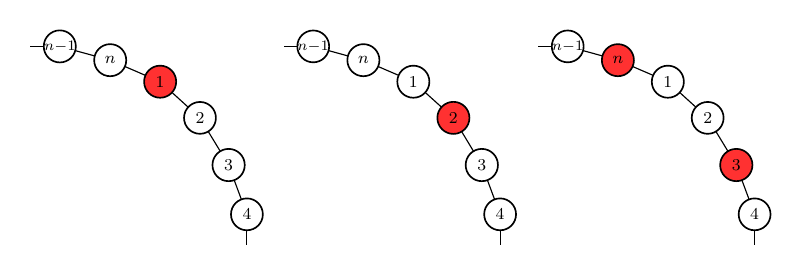}
\caption{Partition of the jamming limits of a cycle.}
\label{fig:particion}
\end{figure}

For the static model, conditioning on the three classes of the partition of $M(C_n)$ yields the convex combination
\begin{equation}
\mathbb{P}\left[N_{C_n}^{\mathrm{stat}} \le k\right]
= \alpha_n\,\mathbb{P}\left[N_{n-3}^{\mathrm{stat}} \le k-1\right]
+ (1 - \alpha_n)\,\mathbb{P}\left[N_{n-6}^{\mathrm{stat}} \le k-2\right],
\label{eq:ciclo_est}
\end{equation}
where $\alpha_n = \frac{2P_{n-2}}{2P_{n-2} + P_{n-5}} \in (0,1)$, since a path of $m$ vertices has $P_{m+1}$ jamming limits.

For the sequential model the situation is considerably simpler. Conditioning on the first occupied vertex --- which is any vertex of the cycle with probability $1/n$ --- that vertex and its two immediate neighbors become blocked, and by symmetry the remaining process develops on a path of $n-3$ vertices. Hence
\begin{equation}
\mathbb{P}\left[N_{C_n}^{\mathrm{seq}} \le k\right]
= \mathbb{P}\left[N_{n-3}^{\mathrm{seq}} \le k-1\right].
\label{eq:ciclo_seq}
\end{equation}

By Lemma~\ref{lem:monotonia_estatica}\textup{(ii)} (with $m = n-3$ and $r = k-1$ when $k \ge 1$; for $k = 0$ both sides vanish), we have $\mathbb{P}\left[N_{n-3}^{\mathrm{stat}} \le k-1\right] \le \mathbb{P}\left[N_{n-6}^{\mathrm{stat}} \le k-2\right]$, and by Theorem~\ref{thm:main_corrected} we also have $\mathbb{P}\left[N_{n-3}^{\mathrm{seq}} \le k-1\right] \le \mathbb{P}\left[N_{n-3}^{\mathrm{stat}} \le k-1\right]$. Combining these two inequalities with \eqref{eq:ciclo_est} and \eqref{eq:ciclo_seq} gives
\[
\begin{aligned}
\mathbb{P}\left[N_{C_n}^{\mathrm{stat}} \le k\right]
&= \alpha_n\,\mathbb{P}\left[N_{n-3}^{\mathrm{stat}} \le k-1\right]
+ (1 - \alpha_n)\,\mathbb{P}\left[N_{n-6}^{\mathrm{stat}} \le k-2\right] \\
&\ge \mathbb{P}\left[N_{n-3}^{\mathrm{stat}} \le k-1\right] \\
&\ge \mathbb{P}\left[N_{n-3}^{\mathrm{seq}} \le k-1\right]
= \mathbb{P}\left[N_{C_n}^{\mathrm{seq}} \le k\right].
\end{aligned}
\]
Therefore, the sequential model dominates the static model on $C_n$.
\end{proof}

\subsection{Computational verification for small graphs}
\label{sec:computacional}
To probe the limits of the property beyond paths and cycles, we computed the laws of $N_G^{\mathrm{seq}}$ and $N_G^{\mathrm{stat}}$ for all small graphs (code in \cite{LopezRivera2026Repositorio}). The algorithm proceeds in three stages. First, a backtracking search with a fixed vertex labeling enumerates all valid selection sequences until reaching a jamming limit. Second, each sequence is mapped to its resulting jamming limit, encoded as a binary vector in $\{0,1\}^n$; the set of distinct vectors determines $M(G)$ and the static law \eqref{eq:est_cdf}. Third, the sequential probability of each jamming limit is computed via \eqref{eq:seq} by summing the path probabilities over all sequences ending in that limit, and grouping by cardinality yields the distribution of $N_G^{\mathrm{seq}}$ in \eqref{eq:seq_cdf}.

\begin{prop}\label{thm:pequenos}
For all connected simple graphs $G = (V,E)$ with $|V| < 6$, the sequential model dominates the static model ($N_G^{\mathrm{stat}} \le_{\mathrm{st}} N_G^{\mathrm{seq}}$).
\end{prop}

\begin{proof}
According to the On-Line Encyclopedia of Integer Sequences (OEIS, sequence A001349 \cite{OEIS_A001349}), there are exactly 31 non-isomorphic connected simple graphs with fewer than 6 vertices. The proof relies on the exact computational enumeration of all configuration spaces and selection probabilities for these 31 graphs. Table~\ref{tab:subset_grafos_acumulada} presents the cumulative distributions obtained for the graphs with richer structures in this family, excluding complete graphs, paths, and cycles (only the values of $k$ in the support of $N_G$ are displayed, as the cumulative distribution functions remain constant at intermediate integers). For a complete verification of all 31 graphs, see the repository \cite{LopezRivera2026Repositorio}.
\end{proof}

\begingroup
\renewcommand{\arraystretch}{1.5}
\begin{longtable}{>{\centering\arraybackslash}m{2.2cm} c c c}
\caption{Cumulative distributions for connected simple graphs with $|V| < 6$ (excluding complete graphs, paths, and cycles).}\label{tab:subset_grafos_acumulada}\\
\toprule
\textbf{Graph} & \textbf{Values of $k$} & \textbf{Sequential distribution} & \textbf{Static distribution}\\
\midrule
\endfirsthead
\toprule
\textbf{Graph} & \textbf{Values of $k$} & \textbf{Sequential distribution} & \textbf{Static distribution}\\
\midrule
\endhead
\bottomrule
\endfoot
\bottomrule
\endlastfoot
\includegraphics[width=1.8cm]{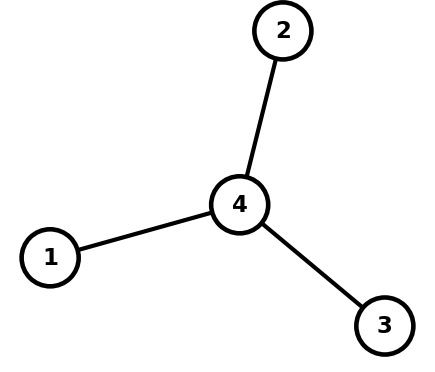} & \begin{tabular}[c]{c} $k=1$ \\ $k=3$ \end{tabular} & \begin{tabular}[c]{c} $\mathbb{P}[N_{G}^{\mathrm{seq}} \le 1] = \frac{1}{4}$ \\ $\mathbb{P}[N_{G}^{\mathrm{seq}} \le 3] = 1$ \end{tabular} & \begin{tabular}[c]{c} $\mathbb{P}[N_{G}^{\mathrm{stat}} \le 1] = \frac{1}{2}$ \\ $\mathbb{P}[N_{G}^{\mathrm{stat}} \le 3] = 1$ \end{tabular}\\
\midrule
\includegraphics[width=1.8cm]{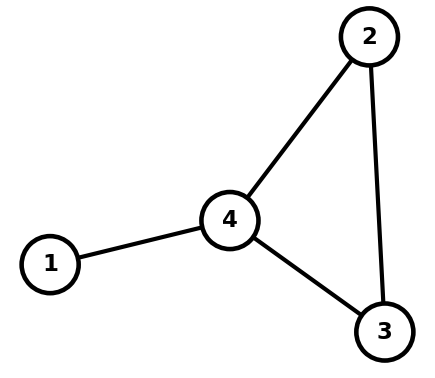} & \begin{tabular}[c]{c} $k=1$ \\ $k=2$ \end{tabular} & \begin{tabular}[c]{c} $\mathbb{P}[N_{G}^{\mathrm{seq}} \le 1] = \frac{1}{4}$ \\ $\mathbb{P}[N_{G}^{\mathrm{seq}} \le 2] = 1$ \end{tabular} & \begin{tabular}[c]{c} $\mathbb{P}[N_{G}^{\mathrm{stat}} \le 1] = \frac{1}{3}$ \\ $\mathbb{P}[N_{G}^{\mathrm{stat}} \le 2] = 1$ \end{tabular}\\
\midrule
\includegraphics[width=1.8cm]{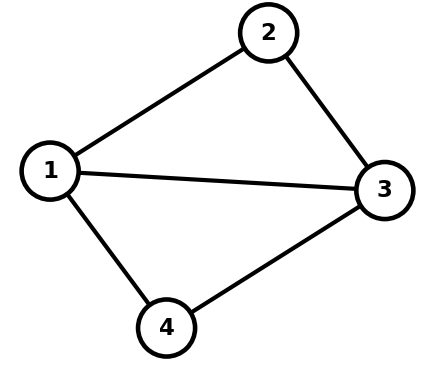} & \begin{tabular}[c]{c} $k=1$ \\ $k=2$ \end{tabular} & \begin{tabular}[c]{c} $\mathbb{P}[N_{G}^{\mathrm{seq}} \le 1] = \frac{1}{2}$ \\ $\mathbb{P}[N_{G}^{\mathrm{seq}} \le 2] = 1$ \end{tabular} & \begin{tabular}[c]{c} $\mathbb{P}[N_{G}^{\mathrm{stat}} \le 1] = \frac{2}{3}$ \\ $\mathbb{P}[N_{G}^{\mathrm{stat}} \le 2] = 1$ \end{tabular}\\
\midrule
\includegraphics[width=1.8cm]{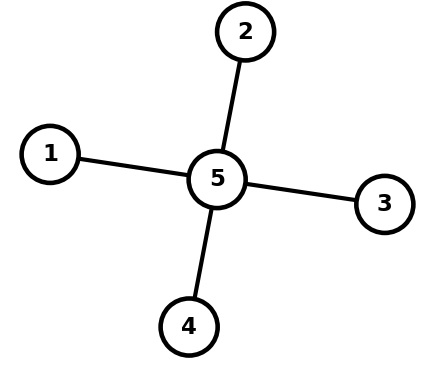} & \begin{tabular}[c]{c} $k=1$ \\ $k=4$ \end{tabular} & \begin{tabular}[c]{c} $\mathbb{P}[N_{G}^{\mathrm{seq}} \le 1] = \frac{1}{5}$ \\ $\mathbb{P}[N_{G}^{\mathrm{seq}} \le 4] = 1$ \end{tabular} & \begin{tabular}[c]{c} $\mathbb{P}[N_{G}^{\mathrm{stat}} \le 1] = \frac{1}{2}$ \\ $\mathbb{P}[N_{G}^{\mathrm{stat}} \le 4] = 1$ \end{tabular}\\
\midrule
\includegraphics[width=1.8cm]{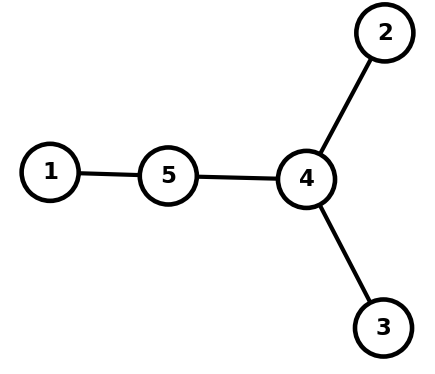} & \begin{tabular}[c]{c} $k=2$ \\ $k=3$ \end{tabular} & \begin{tabular}[c]{c} $\mathbb{P}[N_{G}^{\mathrm{seq}} \le 2] = \frac{4}{15}$ \\ $\mathbb{P}[N_{G}^{\mathrm{seq}} \le 3] = 1$ \end{tabular} & \begin{tabular}[c]{c} $\mathbb{P}[N_{G}^{\mathrm{stat}} \le 2] = \frac{1}{3}$ \\ $\mathbb{P}[N_{G}^{\mathrm{stat}} \le 3] = 1$ \end{tabular}\\
\midrule
\includegraphics[width=1.8cm]{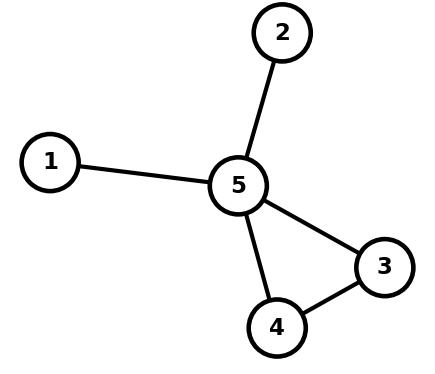} & \begin{tabular}[c]{c} $k=1$ \\ $k=3$ \end{tabular} & \begin{tabular}[c]{c} $\mathbb{P}[N_{G}^{\mathrm{seq}} \le 1] = \frac{1}{5}$ \\ $\mathbb{P}[N_{G}^{\mathrm{seq}} \le 3] = 1$ \end{tabular} & \begin{tabular}[c]{c} $\mathbb{P}[N_{G}^{\mathrm{stat}} \le 1] = \frac{1}{3}$ \\ $\mathbb{P}[N_{G}^{\mathrm{stat}} \le 3] = 1$ \end{tabular}\\
\midrule
\includegraphics[width=1.8cm]{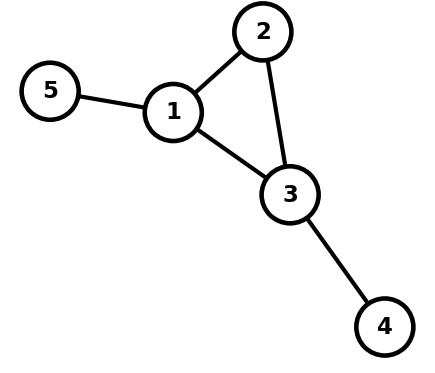} & \begin{tabular}[c]{c} $k=2$ \\ $k=3$ \end{tabular} & \begin{tabular}[c]{c} $\mathbb{P}[N_{G}^{\mathrm{seq}} \le 2] = \frac{8}{15}$ \\ $\mathbb{P}[N_{G}^{\mathrm{seq}} \le 3] = 1$ \end{tabular} & \begin{tabular}[c]{c} $\mathbb{P}[N_{G}^{\mathrm{stat}} \le 2] = \frac{2}{3}$ \\ $\mathbb{P}[N_{G}^{\mathrm{stat}} \le 3] = 1$ \end{tabular}\\
\midrule
\includegraphics[width=1.8cm]{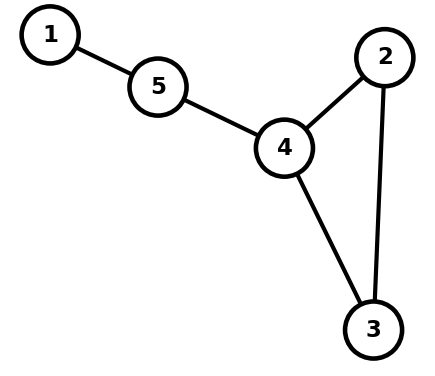} & $k=2$ & $\mathbb{P}[N_{G}^{\mathrm{seq}} \le 2] = 1$ & $\mathbb{P}[N_{G}^{\mathrm{stat}} \le 2] = 1$\\
\midrule
\includegraphics[width=1.8cm]{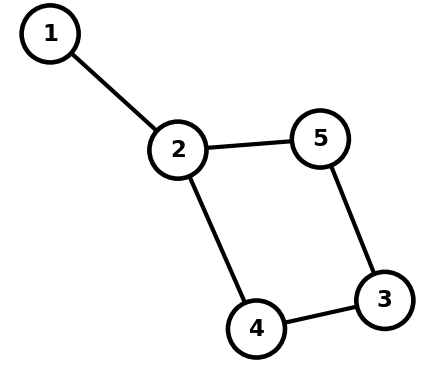} & \begin{tabular}[c]{c} $k=2$ \\ $k=3$ \end{tabular} & \begin{tabular}[c]{c} $\mathbb{P}[N_{G}^{\mathrm{seq}} \le 2] = \frac{7}{15}$ \\ $\mathbb{P}[N_{G}^{\mathrm{seq}} \le 3] = 1$ \end{tabular} & \begin{tabular}[c]{c} $\mathbb{P}[N_{G}^{\mathrm{stat}} \le 2] = \frac{2}{3}$ \\ $\mathbb{P}[N_{G}^{\mathrm{stat}} \le 3] = 1$ \end{tabular}\\
\midrule
\includegraphics[width=1.8cm]{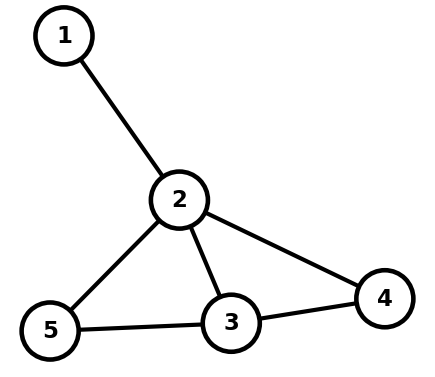} & \begin{tabular}[c]{c} $k=1$ \\ $k=2$ \\ $k=3$ \end{tabular} & \begin{tabular}[c]{c} $\mathbb{P}[N_{G}^{\mathrm{seq}} \le 1] = \frac{1}{5}$ \\ $\mathbb{P}[N_{G}^{\mathrm{seq}} \le 2] = \frac{7}{15}$ \\ $\mathbb{P}[N_{G}^{\mathrm{seq}} \le 3] = 1$ \end{tabular} & \begin{tabular}[c]{c} $\mathbb{P}[N_{G}^{\mathrm{stat}} \le 1] = \frac{1}{3}$ \\ $\mathbb{P}[N_{G}^{\mathrm{stat}} \le 2] = \frac{2}{3}$ \\ $\mathbb{P}[N_{G}^{\mathrm{stat}} \le 3] = 1$ \end{tabular}\\
\midrule
\includegraphics[width=1.8cm]{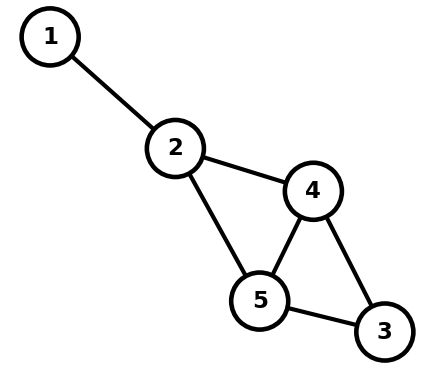} & $k=2$ & $\mathbb{P}[N_{G}^{\mathrm{seq}} \le 2] = 1$ & $\mathbb{P}[N_{G}^{\mathrm{stat}} \le 2] = 1$\\
\midrule
\includegraphics[width=1.8cm]{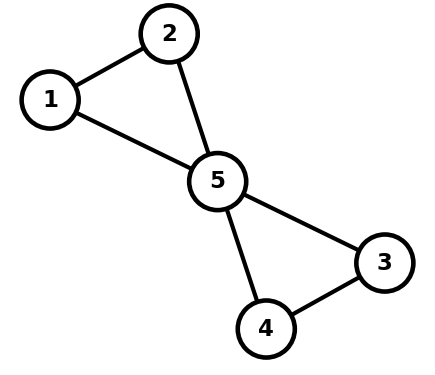} & \begin{tabular}[c]{c} $k=1$ \\ $k=2$ \end{tabular} & \begin{tabular}[c]{c} $\mathbb{P}[N_{G}^{\mathrm{seq}} \le 1] = \frac{1}{5}$ \\ $\mathbb{P}[N_{G}^{\mathrm{seq}} \le 2] = 1$ \end{tabular} & \begin{tabular}[c]{c} $\mathbb{P}[N_{G}^{\mathrm{stat}} \le 1] = \frac{1}{5}$ \\ $\mathbb{P}[N_{G}^{\mathrm{stat}} \le 2] = 1$ \end{tabular}\\
\midrule
\includegraphics[width=1.8cm]{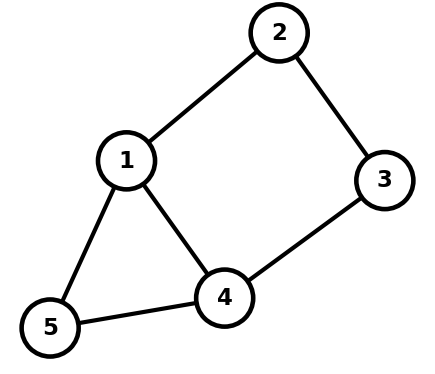} & $k=2$ & $\mathbb{P}[N_{G}^{\mathrm{seq}} \le 2] = 1$ & $\mathbb{P}[N_{G}^{\mathrm{stat}} \le 2] = 1$\\
\midrule
\includegraphics[width=1.8cm]{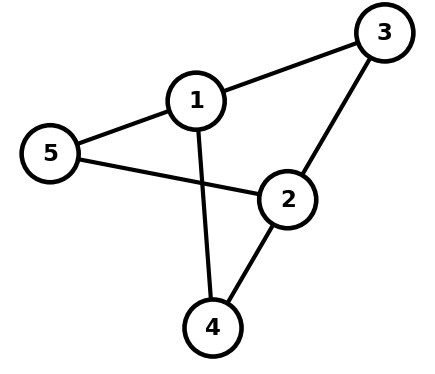} & \begin{tabular}[c]{c} $k=2$ \\ $k=3$ \end{tabular} & \begin{tabular}[c]{c} $\mathbb{P}[N_{G}^{\mathrm{seq}} \le 2] = \frac{2}{5}$ \\ $\mathbb{P}[N_{G}^{\mathrm{seq}} \le 3] = 1$ \end{tabular} & \begin{tabular}[c]{c} $\mathbb{P}[N_{G}^{\mathrm{stat}} \le 2] = \frac{1}{2}$ \\ $\mathbb{P}[N_{G}^{\mathrm{stat}} \le 3] = 1$ \end{tabular}\\
\midrule
\includegraphics[width=1.8cm]{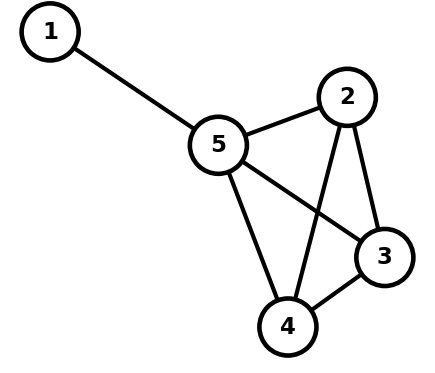} & \begin{tabular}[c]{c} $k=1$ \\ $k=2$ \end{tabular} & \begin{tabular}[c]{c} $\mathbb{P}[N_{G}^{\mathrm{seq}} \le 1] = \frac{1}{5}$ \\ $\mathbb{P}[N_{G}^{\mathrm{seq}} \le 2] = 1$ \end{tabular} & \begin{tabular}[c]{c} $\mathbb{P}[N_{G}^{\mathrm{stat}} \le 1] = \frac{1}{4}$ \\ $\mathbb{P}[N_{G}^{\mathrm{stat}} \le 2] = 1$ \end{tabular}\\
\midrule
\includegraphics[width=1.8cm]{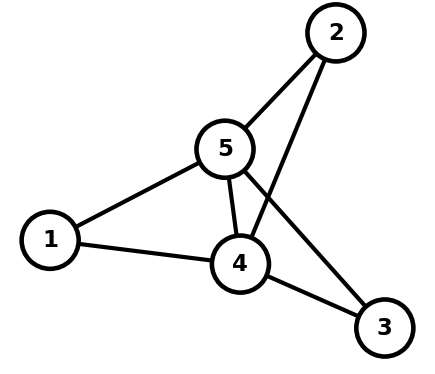} & \begin{tabular}[c]{c} $k=1$ \\ $k=3$ \end{tabular} & \begin{tabular}[c]{c} $\mathbb{P}[N_{G}^{\mathrm{seq}} \le 1] = \frac{2}{5}$ \\ $\mathbb{P}[N_{G}^{\mathrm{seq}} \le 3] = 1$ \end{tabular} & \begin{tabular}[c]{c} $\mathbb{P}[N_{G}^{\mathrm{stat}} \le 1] = \frac{2}{3}$ \\ $\mathbb{P}[N_{G}^{\mathrm{stat}} \le 3] = 1$ \end{tabular}\\
\midrule
\includegraphics[width=1.8cm]{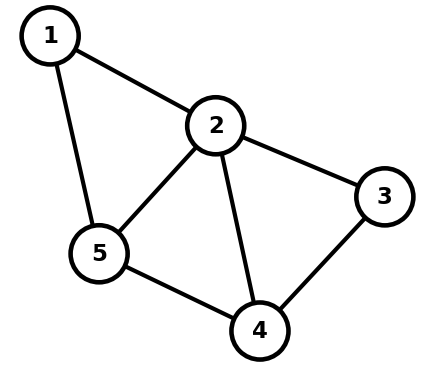} & \begin{tabular}[c]{c} $k=1$ \\ $k=2$ \end{tabular} & \begin{tabular}[c]{c} $\mathbb{P}[N_{G}^{\mathrm{seq}} \le 1] = \frac{1}{5}$ \\ $\mathbb{P}[N_{G}^{\mathrm{seq}} \le 2] = 1$ \end{tabular} & \begin{tabular}[c]{c} $\mathbb{P}[N_{G}^{\mathrm{stat}} \le 1] = \frac{1}{4}$ \\ $\mathbb{P}[N_{G}^{\mathrm{stat}} \le 2] = 1$ \end{tabular}\\
\midrule
\includegraphics[width=1.8cm]{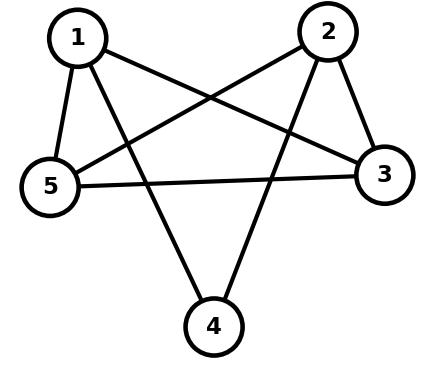} & $k=2$ & $\mathbb{P}[N_{G}^{\mathrm{seq}} \le 2] = 1$ & $\mathbb{P}[N_{G}^{\mathrm{stat}} \le 2] = 1$\\
\midrule
\includegraphics[width=1.8cm]{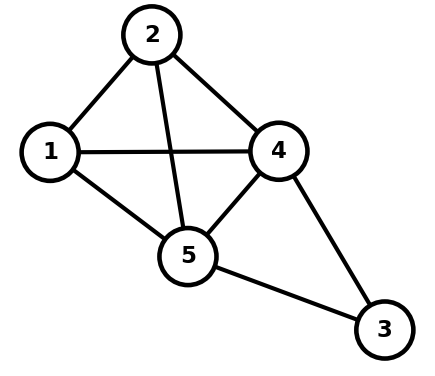} & \begin{tabular}[c]{c} $k=1$ \\ $k=2$ \end{tabular} & \begin{tabular}[c]{c} $\mathbb{P}[N_{G}^{\mathrm{seq}} \le 1] = \frac{2}{5}$ \\ $\mathbb{P}[N_{G}^{\mathrm{seq}} \le 2] = 1$ \end{tabular} & \begin{tabular}[c]{c} $\mathbb{P}[N_{G}^{\mathrm{stat}} \le 1] = \frac{1}{2}$ \\ $\mathbb{P}[N_{G}^{\mathrm{stat}} \le 2] = 1$ \end{tabular}\\
\midrule
\includegraphics[width=1.8cm]{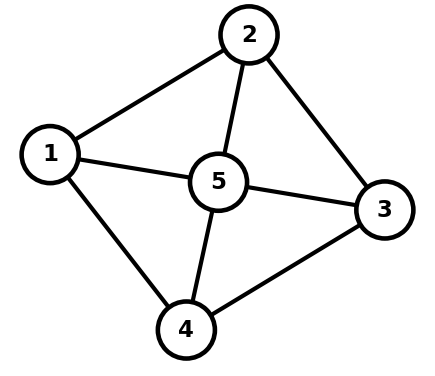} & \begin{tabular}[c]{c} $k=1$ \\ $k=2$ \end{tabular} & \begin{tabular}[c]{c} $\mathbb{P}[N_{G}^{\mathrm{seq}} \le 1] = \frac{1}{5}$ \\ $\mathbb{P}[N_{G}^{\mathrm{seq}} \le 2] = 1$ \end{tabular} & \begin{tabular}[c]{c} $\mathbb{P}[N_{G}^{\mathrm{stat}} \le 1] = \frac{1}{3}$ \\ $\mathbb{P}[N_{G}^{\mathrm{stat}} \le 2] = 1$ \end{tabular}\\
\midrule
\includegraphics[width=1.8cm]{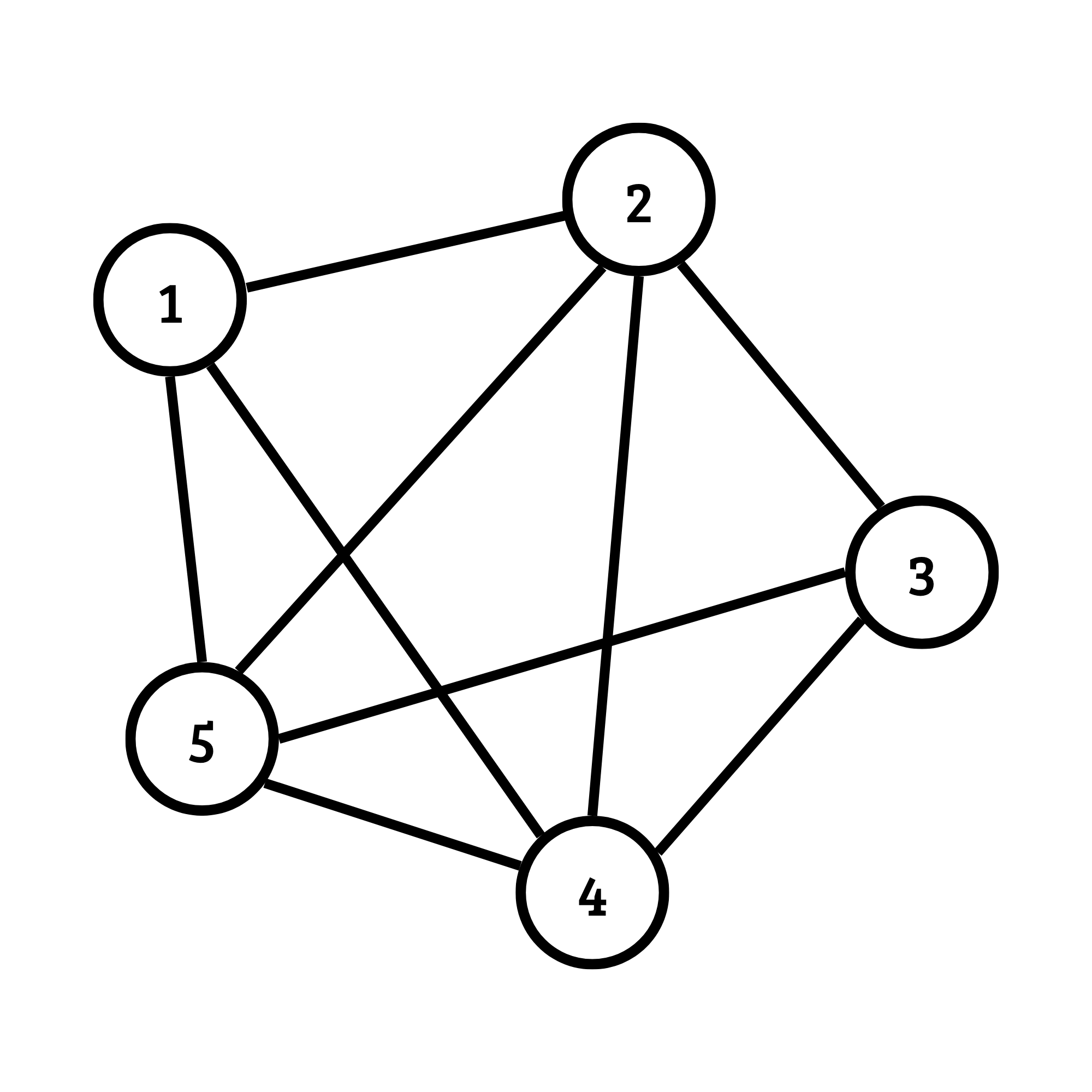} & \begin{tabular}[c]{c} $k=1$ \\ $k=2$ \end{tabular} & \begin{tabular}[c]{c} $\mathbb{P}[N_{G}^{\mathrm{seq}} \le 1] = \frac{3}{5}$ \\ $\mathbb{P}[N_{G}^{\mathrm{seq}} \le 2] = 1$ \end{tabular} & \begin{tabular}[c]{c} $\mathbb{P}[N_{G}^{\mathrm{stat}} \le 1] = \frac{3}{4}$ \\ $\mathbb{P}[N_{G}^{\mathrm{stat}} \le 2] = 1$ \end{tabular}\\
\bottomrule
\end{longtable}
\endgroup

\section{Counterexamples and \texorpdfstring{\textcolor{revcolor}{Open Questions}}{Open Questions}}
\label{sec:counterexamples}

While the results of Sections~\ref{sec:paths} and~\ref{sec:cycles_small} establish $N_G^{\mathrm{stat}} \le_{\mathrm{st}} N_G^{\mathrm{seq}}$ for paths, cycles, and connected graphs with fewer than 6 vertices, this dominance does not hold on arbitrary graphs. Below we construct explicit counterexamples with $|V| \ge 6$ and $\Delta(G) \ge 4$ on which the static model strictly dominates the sequential one, show that while stochastic comparability holds for all connected graphs with up to 7 vertices it fails starting at 8 vertices, \rev{present counterexamples with maximum degree 3, and state some open questions}.

\subsection{A family of counterexamples}
\label{sec:familia}
\begin{figure}[ht]
\centering
\includegraphics[width=0.5\textwidth]{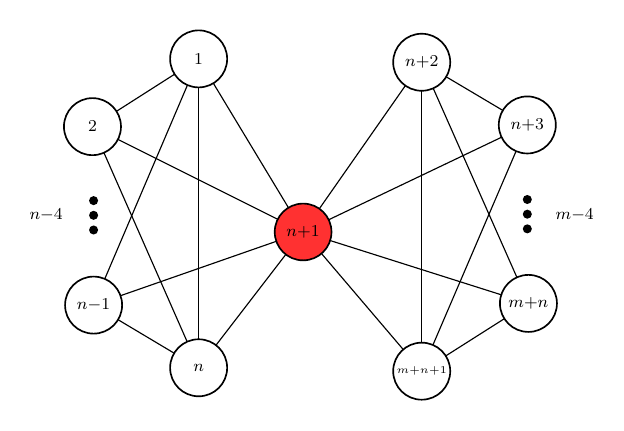}
\caption{Family of graphs where the sequential model does not dominate the static model.}
\label{fig:contraejemplo_general}
\end{figure}

Starting from 6 vertices, it is always possible to find a graph on which the sequential model does not dominate the static model. Consider the graph $G$ formed by two disjoint complete graphs, one with $n$ vertices and one with $m$ vertices, together with a central vertex adjacent to all vertices of both complete graphs (Figure~\ref{fig:contraejemplo_general}); throughout this subsection, $n$ and $m$ denote the sizes of the two complete graphs, so that $|V| = n+m+1$. By construction, no vertex of the first complete graph is connected to any vertex of the second.

Due to this structure, the only possible jamming limits are those with a single occupied seat --- the central vertex --- and those with exactly two occupied seats, one belonging to each complete graph. Hence $N_G^{\mathrm{seq}}$ and $N_G^{\mathrm{stat}}$ only take the values $1$ and $2$, and $\mathbb{P}\left[N_G^{\mathrm{stat}} \le 2\right] = \mathbb{P}\left[N_G^{\mathrm{seq}} \le 2\right] = 1$. To show that $N_G^{\mathrm{seq}} \le_{\mathrm{st}} N_G^{\mathrm{stat}}$ strictly, it suffices to verify that $\mathbb{P}\left[N_G^{\mathrm{seq}} = 1\right] > \mathbb{P}\left[N_G^{\mathrm{stat}} = 1\right]$.

In the sequential model, the process ends with a single occupied seat exactly when the first selected vertex is the central one, among $n+m+1$ available vertices; hence
\[
\mathbb{P}\left[N_G^{\mathrm{seq}} = 1\right] = \frac{1}{n+m+1}.
\]
In the static model there is exactly one jamming limit with a single occupied seat, corresponding to the central vertex, while a jamming limit with two occupied seats is obtained by choosing one vertex from each complete graph, giving $nm$ such limits. Hence $|M(G)| = nm+1$ and
\[
\mathbb{P}\left[N_G^{\mathrm{stat}} = 1\right] = \frac{1}{nm+1}.
\]

Therefore, the static model strictly dominates the sequential model whenever
\[
\frac{1}{n+m+1} > \frac{1}{nm+1} \iff (n-1)(m-1) > 1.
\]
Fixing $n = 2$, the condition reduces to $m > 2$. Thus, for every $m \ge 3$, the static model strictly dominates the sequential one ($N_G^{\mathrm{seq}} \le_{\mathrm{st}} N_G^{\mathrm{stat}}$) on this graph, which has $|V| = m+3 \ge 6$ vertices and maximum degree $\Delta(G) = n+m \ge 5$.

\subsection{A counterexample with maximum degree 4}
\label{sec:grado4}
\begin{figure}[ht]
\centering
\includegraphics[width=0.3\textwidth]{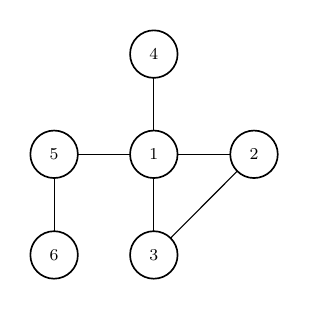}
\caption{Graph of maximum degree 4 where the sequential model does not dominate the static model.}
\label{fig:contraejemplo_particular}
\end{figure}

\begin{figure}[ht]
\centering
\includegraphics[width=0.7\textwidth]{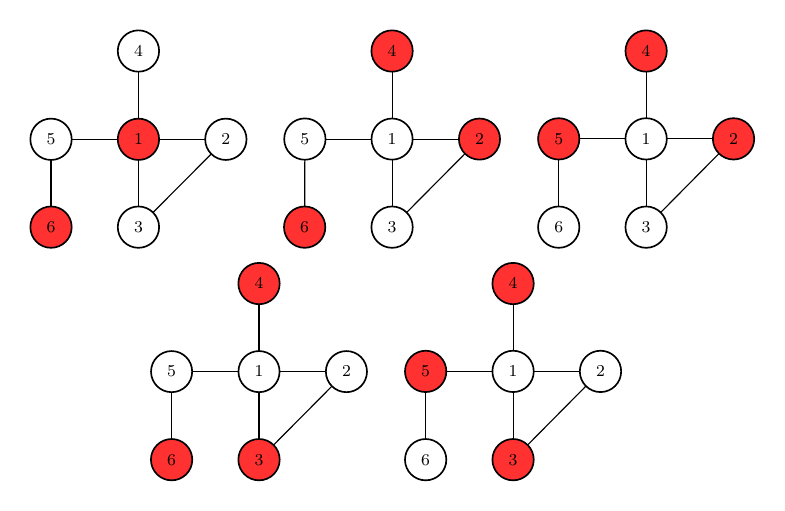}
\caption{Jamming limits of the graph of Figure~\ref{fig:contraejemplo_particular}.}
\label{fig:contraejemplo_particular_1}
\end{figure}

The family of counterexamples presented above shows that for every $\Delta(G) \ge 5$ there exists a graph on which the static model strictly dominates the sequential model. It remains to analyze the case $\Delta(G) = 4$. Consider the graph $G$ with vertices $V = \{1, 2, 3, 4, 5, 6\}$ and edges arranged so that vertex $1$ is connected to vertices $2, 3, 4$, and $5$; vertex $5$ is connected to $6$; and vertex $2$ is connected to $3$ (Figure~\ref{fig:contraejemplo_particular}). Its jamming limits are shown in Figure~\ref{fig:contraejemplo_particular_1}: there are exactly five, namely $\{1,6\}$ of size 2, and $\{2,4,5\}$, $\{2,4,6\}$, $\{3,4,5\}$, $\{3,4,6\}$ of size 3. Hence $N_G^{\mathrm{seq}}$ and $N_G^{\mathrm{stat}}$ take only the values $2$ and $3$.

For the process to end with exactly two occupied vertices, these must be vertices $1$ and $6$. If the first selected vertex is $6$, then vertices $5$ and $6$ are no longer available, leaving 4 vertices that can be chosen in the next step, one of which is vertex $1$. Hence the probability of the sequence $(6,1)$ is
\[
\frac{1}{6} \cdot \frac{1}{4} = \frac{1}{24}.
\]
If the first selected vertex is $1$, then all other vertices are blocked except $6$, and the probability of the sequence $(1,6)$ is $1/6$. Therefore,
\[
\mathbb{P}\left[N_G^{\mathrm{seq}} \le 2\right] = \mathbb{P}\left[N_G^{\mathrm{seq}} = 2\right] = \frac{1}{6} + \frac{1}{24} = \frac{5}{24}.
\]
In the static model there are five equiprobable jamming limits and only one of them contains two occupied vertices, so
\[
\mathbb{P}\left[N_G^{\mathrm{stat}} \le 2\right] = \mathbb{P}\left[N_G^{\mathrm{stat}} = 2\right] = \frac{1}{5}.
\]
Since $5/24 > 1/5$, the static model strictly dominates the sequential model ($N_G^{\mathrm{seq}} \le_{\mathrm{st}} N_G^{\mathrm{stat}}$) on this graph of maximum degree 4.

\subsection{Stochastic comparability fails}
\label{sec:comparabilidad}
So far, in all graphs presented as counterexamples to $N_G^{\mathrm{stat}} \le_{\mathrm{st}} N_G^{\mathrm{seq}}$, the random variable $N_G$ takes only two possible values $k_1 < k_2$. Whenever the support of $N_G$ consists of at most two values, the cumulative distribution functions have only one nontrivial value before reaching $1$ (namely, at $k = k_1$), which forces one of the two models to stochastically dominate the other. Consequently, a necessary condition for the two models to be stochastically incomparable is that $G$ admits maximal independent sets of at least three different cardinalities. This leads naturally to asking whether on every graph at least one of the two models dominates the other.

The answer is negative. On all connected graphs with up to 5 vertices, the sequential model dominates the static model (Proposition~\ref{thm:pequenos}), and our computational verification (code in \cite{LopezRivera2026Repositorio}) establishes the following observation for graphs with up to 7 vertices.

\begin{observation}\label{obs:7vertices}
On all connected simple graphs with up to 7 vertices --- comprising 112 non-isomorphic connected graphs with 6 vertices and 853 with 7 vertices (sequence A001349 of the OEIS \cite{OEIS_A001349}) --- one of the two models first-order stochastically dominates the other.
\end{observation}

More precisely, among the 112 connected graphs with 6 vertices, the two laws coincide on 31 of them, the sequential model strictly dominates on 79, and the static model strictly dominates on 2; among the 853 connected graphs with 7 vertices, the corresponding numbers are 119, 701 and 33. Every graph in these two families on which the static model strictly dominates has maximum degree at least 4.

\begin{figure}[ht]
\centering
\includegraphics[width=0.5\textwidth]{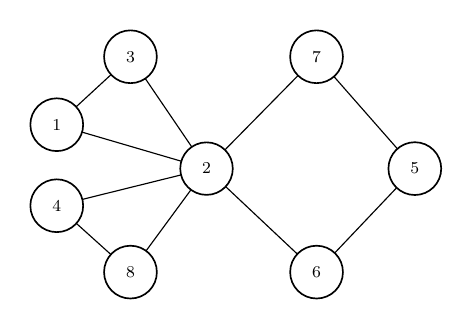}
\caption{An 8-vertex graph on which neither model dominates the other.}
\label{fig:pez}
\end{figure}

For 8 vertices, however, universal comparability fails. An example of an 8-vertex connected graph on which neither model dominates the other is the graph $G$ of Figure~\ref{fig:pez}, formed by two triangles, $\{1,2,3\}$ and $\{2,4,8\}$, and a $4$-cycle, $(2,7,5,6)$, sharing the central vertex $2$. This graph has $9$ jamming limits: one of size $2$ (namely $\{2,5\}$), four of size $3$ (of the form $\{u,v,5\}$ with $u \in \{1,3\}$ and $v \in \{4,8\}$), and four of size $4$ (of the form $\{u,v,6,7\}$ with $u \in \{1,3\}$ and $v \in \{4,8\}$). As shown in Table~\ref{tab:pez}, at $k=2$ we have $\mathbb{P}[N_G^{\mathrm{stat}} \le 2] = 1/9 < 3/20 = \mathbb{P}[N_G^{\mathrm{seq}} \le 2]$, whereas at $k=3$ the inequality reverses to $\mathbb{P}[N_G^{\mathrm{stat}} \le 3] = 5/9 > 5/12 = \mathbb{P}[N_G^{\mathrm{seq}} \le 3]$. Since the cumulative distribution functions cross, neither model dominates the other.

\begin{table}[H]
\centering
\renewcommand{\arraystretch}{1.25}
\begin{tabular}{ccc}
\toprule
$k$ & $\mathbb{P}[N_G^{\mathrm{stat}}\leq k]$ & $\mathbb{P}[N_G^{\mathrm{seq}}\leq k]$\\
\midrule
$2$ & $1/9$ & $3/20$\\
$3$ & $5/9$ & $5/12$\\
$4$ & $1$ & $1$\\
\bottomrule
\end{tabular}
\caption{Cumulative distributions of the number of occupied vertices for the 8-vertex graph of Figure~\ref{fig:pez}.}
\label{tab:pez}
\end{table}

\subsection{\texorpdfstring{\textcolor{revcolor}{Counterexamples with maximum degree 3}}{Counterexamples with maximum degree 3}}
\label{sec:grado3}
\rev{Every graph presented so far on which the sequential model fails to dominate the static one has a vertex of degree at least 4. One might therefore expect the sequential model to dominate the static one on every graph with maximum degree at most 3. To see that this is not the case, we consider the following two examples, which were provided to us by D.~Jiang (personal communication, 2026). The second one, a tree, was also found independently by A.~Jangid (personal communication, 2026).}

\rev{Let $G$ be the connected unicyclic graph with 15 vertices and maximum degree 3 shown in Figure~\ref{fig:jiangG}. It has 55 jamming limits, one of size 6 and 54 of size 7, so that the number of occupied vertices takes only the values 6 and 7; its laws in the two models are given in Table~\ref{tab:jiangG}. Since
\[
\mathbb{P}\left[N_G^{\mathrm{seq}} \le 6\right] = \frac{64}{3375} \approx 0.018963 > 0.018182 \approx \frac{1}{55} = \mathbb{P}\left[N_G^{\mathrm{stat}} \le 6\right],
\]
the static model strictly dominates the sequential model ($N_G^{\mathrm{seq}} \le_{\mathrm{st}} N_G^{\mathrm{stat}}$) on this graph of maximum degree 3. In particular, $\mathbb{E}[N_G^{\mathrm{seq}}] - \mathbb{E}[N_G^{\mathrm{stat}}] = -29/37125 < 0$, so that even the inequality between the mean occupancies fails on $G$. The graph obtained from $G$ by replacing the edge $\{6,14\}$ with $\{8,14\}$ has the same two laws.}

\begin{figure}[!ht]
\centering
\begin{tikzpicture}[scale=0.85,
  vx/.style={circle, draw, inner sep=0pt, minimum size=5.5mm, font=\scriptsize},
  occ/.style={vx, fill=gray!40}]
\node[occ] (n4)  at ( 0.73, 2.70) {$4$};
\node[vx]  (n7)  at (-0.75, 1.83) {$7$};
\node[vx]  (n11) at ( 0.73, 1.83) {$11$};
\node[occ] (n0)  at ( 0.00, 1.00) {$0$};
\node[vx]  (n14) at ( 0.00, 0.00) {$14$};
\node[vx]  (n6)  at (-0.95,-0.68) {$6$};
\node[vx]  (n1)  at ( 0.93,-0.68) {$1$};
\node[occ] (n12) at ( 1.95,-0.38) {$12$};
\node[vx]  (n2)  at ( 2.83, 0.13) {$2$};
\node[occ] (n9)  at ( 3.70, 0.63) {$9$};
\node[vx]  (n5)  at ( 2.83,-0.83) {$5$};
\node[occ] (n13) at (-0.60,-1.75) {$13$};
\node[vx]  (n8)  at ( 0.58,-1.75) {$8$};
\node[vx]  (n3)  at (-1.20,-2.58) {$3$};
\node[occ] (n10) at (-1.83,-3.40) {$10$};
\draw (n7)--(n0)--(n11)--(n4) (n0)--(n14);
\draw (n14)--(n6)--(n13)--(n8)--(n1)--(n14);
\draw (n1)--(n12)--(n2)--(n9) (n12)--(n5);
\draw (n13)--(n3)--(n10);
\end{tikzpicture}
\caption{\textcolor{revcolor}{The unicyclic graph $G$, with 15 vertices and maximum degree 3. Shaded: its unique jamming limit of size 6, $\{0, 4, 9, 10, 12, 13\}$.}}
\label{fig:jiangG}

\bigskip
{\color{revcolor}\renewcommand{\arraystretch}{1.25}
\begin{tabular}{ccc}
\toprule
$k$ & $\mathbb{P}[N_G^{\mathrm{seq}} = k]$ & $\mathbb{P}[N_G^{\mathrm{stat}} = k]$\\
\midrule
$6$ & $64/3375$ & $1/55$\\
$7$ & $3311/3375$ & $54/55$\\
\bottomrule
\end{tabular}}
\captionof{table}{\textcolor{revcolor}{Laws of the number of occupied vertices for the graph $G$ of Figure~\ref{fig:jiangG}.}}
\label{tab:jiangG}
\vspace{4pt}
\end{figure}

\begin{figure}[!ht]
\centering
\begin{tikzpicture}[scale=0.9,
  vx/.style={circle, draw, inner sep=0pt, minimum size=5.5mm, font=\scriptsize},
  occ/.style={vx, fill=gray!40}]
\foreach \lab/\x/\s in {10/0/vx, 7/1/occ, 4/2/vx, 1/3/vx, 0/4/occ, 2/5/vx, 5/6/vx, 8/7/occ, 12/8/vx}
  \node[\s] (t\lab) at (\x, 0) {$\lab$};
\node[vx]  (t9)  at (1,-1) {$9$};  \node[occ] (t13) at (1,-2) {$13$};
\node[vx]  (t3)  at (4,-1) {$3$};  \node[occ] (t6)  at (4,-2) {$6$};
\node[vx]  (t11) at (7,-1) {$11$}; \node[occ] (t14) at (7,-2) {$14$};
\draw (t10)--(t7)--(t4)--(t1)--(t0)--(t2)--(t5)--(t8)--(t12);
\draw (t7)--(t9)--(t13) (t0)--(t3)--(t6) (t8)--(t11)--(t14);
\end{tikzpicture}
\caption{\textcolor{revcolor}{The tree $T$, with 15 vertices and maximum degree 3. Shaded: its unique jamming limit of size 6.}}
\label{fig:jiangT}

\bigskip
{\color{revcolor}\renewcommand{\arraystretch}{1.25}
\begin{tabular}{ccc}
\toprule
$k$ & $\mathbb{P}[N_T^{\mathrm{seq}} = k]$ & $\mathbb{P}[N_T^{\mathrm{stat}} = k]$\\
\midrule
$6$ & $64/3375$ & $1/55$\\
$7$ & $7993019/10135125$ & $10/11$\\
$8$ & $1949914/10135125$ & $4/55$\\
\bottomrule
\end{tabular}}
\captionof{table}{\textcolor{revcolor}{Laws of the number of occupied vertices for the tree $T$ of Figure~\ref{fig:jiangT}.}}
\label{tab:jiangT}
\end{figure}

\rev{Now let $T$ be the tree with 15 vertices and maximum degree 3 shown in Figure~\ref{fig:jiangT}. It also has 55 jamming limits: one of size 6, 50 of size 7 and 4 of size 8 (Table~\ref{tab:jiangT}). In this case,
\begin{gather*}
\mathbb{P}\left[N_T^{\mathrm{seq}} \le 6\right] = \frac{64}{3375} > \frac{1}{55} = \mathbb{P}\left[N_T^{\mathrm{stat}} \le 6\right],\\
\mathbb{P}\left[N_T^{\mathrm{seq}} \ge 8\right] = \frac{1949914}{10135125} > \frac{4}{55} = \mathbb{P}\left[N_T^{\mathrm{stat}} \ge 8\right].
\end{gather*}
The first inequality shows that the sequential model does not dominate the static one on $T$, and the second one that the static model does not dominate the sequential one either. Therefore, a maximum degree of at most 3 guarantees neither the dominance of the sequential model nor the comparability of the two models, not even on trees.}

\rev{According to Jiang, the graphs above are the smallest counterexamples with maximum degree 3. We have verified the two examples in exact rational arithmetic and confirmed their minimality by an exhaustive computation (code in \cite{LopezRivera2026Repositorio}): the sequential model dominates the static one on every connected graph with maximum degree at most 3 and at most 14 vertices ($358{,}891$ non-isomorphic graphs in total), and among the $1{,}016{,}740$ such graphs with 15 vertices the dominance fails on exactly three, namely the graph $G$, its variant with the edge $\{8,14\}$, and the tree $T$. In contrast, Jiang found no violation of the dominance among the connected 3-regular graphs with at most 18 vertices.}

\subsection{\texorpdfstring{\textcolor{revcolor}{Open questions}}{Open questions}}
\label{sec:preguntas}
\rev{In summary, the sequential model dominates the static one on all paths and cycles, on all connected graphs with fewer than 6 vertices, on all connected graphs with maximum degree at most 3 and at most 14 vertices, and, by Corollary~\ref{cor:conexidad}, on all disjoint unions of such graphs. On the other hand, the dominance fails on some graphs with maximum degree at least 4 from 6 vertices on, and on some graphs with maximum degree 3 and 15 vertices, including a tree. Hence neither the number of vertices nor the maximum degree alone determines the direction of the comparison, and we are led to the following questions.}
\begin{itemize}
\item[\textcolor{revcolor}{\textbullet}] \rev{Which structural properties of a finite graph $G$ ensure that the sequential model dominates the static model?}
\item[\textcolor{revcolor}{\textbullet}] \rev{More modestly, on which graphs does the inequality $\mathbb{E}[N_G^{\mathrm{stat}}] \le \mathbb{E}[N_G^{\mathrm{seq}}]$ between the mean occupancies hold? The unicyclic graph $G$ of Figure~\ref{fig:jiangG} shows that it may fail even when the maximum degree is 3.}
\end{itemize}

\begin{figure}[ht]
\centering
\includegraphics[width=0.7\textwidth]{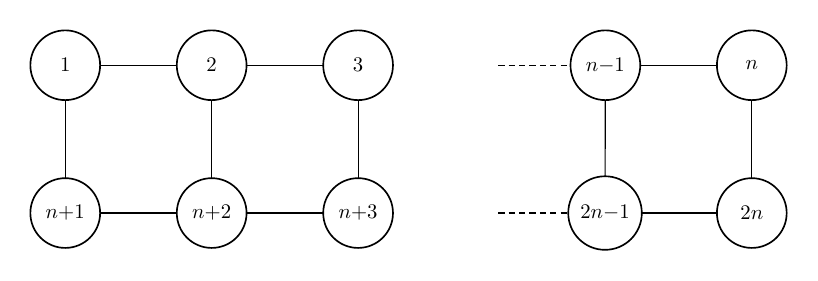}
\caption{The $2 \times n$ grid graph.}
\label{fig:mesa}
\end{figure}

\rev{A natural family on which these questions remain open is that of ladder graphs.} In the dining table process \cite{Dinning-table} (see also \cite{Baram1992, Fan1992, Oh2017} for related studies on ladder and grid graphs), the underlying graphs are the grids $2 \times n$ (Figure~\ref{fig:mesa}). There, the asymptotic density of occupied seats in the sequential model is $(2-e^{-1})/4 \approx 0.4080$, which differs from the value $(1-e^{-2})/2 \approx 0.4323$ of the linear case, whereas the density obtained when all jamming limits are equiprobable is $1/(5-\sqrt{5}) \approx 0.3618$, which is strictly smaller \cite{Dinning-table}.
\rev{These grids have maximum degree at most 3 and, for $n \le 7$, they are among the graphs covered by the exhaustive computation of Section~\ref{sec:grado3}, so the sequential model dominates the static one on them. We do not know whether this is the case for every $n$. If it is, then} for every $n$ the mean occupancy of the sequential model is at least that of the static model, a finite-$n$ counterpart of the asymptotic comparison above.

\section{Disconnected Graphs}
\label{sec:desconexos}
We now analyze how stochastic dominance behaves on a disconnected graph in terms of the occupancy distributions of its connected components.

\subsection{Convolution and inheritance}
\label{sec:convolucion}
\begin{prop}\label{prop:convolucion}
Let $G$ be a finite simple disconnected graph with connected components $G_1, G_2, \dots, G_r$. Then, in both the sequential and the static model, the random variables $N_{G_1}, \dots, N_{G_r}$ are independent and $N_G = \sum_{j=1}^{r} N_{G_j}$. Consequently,
\[
\mathbb{P}\left[N_G = k\right] = \sum_{i_1 + i_2 + \cdots + i_r = k} \prod_{j=1}^{r} \mathbb{P}\left[N_{G_j} = i_j\right]
\]
and
\[
\mathbb{P}\left[N_G \le k\right] = \sum_{m=0}^{k} \sum_{i_1 + i_2 + \cdots + i_r = m} \prod_{j=1}^{r} \mathbb{P}\left[N_{G_j} = i_j\right].
\]
\end{prop}

\begin{proof}
Since there are no edges between distinct connected components, a subset $I \subseteq V(G)$ is a maximal independent set of $G$ if and only if $I = I_1 \cup \dots \cup I_r$, where each $I_j = I \cap V(G_j)$ is a maximal independent set of $G_j$. Thus $M(G)$ is in bijection with $M(G_1) \times \dots \times M(G_r)$ and $|I| = \sum_{j=1}^{r} |I_j|$. In the static model, the uniform measure on $M(G)$ is the product of the uniform measures on $M(G_j)$, so $N_{G_1}^{\mathrm{stat}}, \dots, N_{G_r}^{\mathrm{stat}}$ are independent.

In the sequential model, selecting available vertices uniformly at random at each step is equivalent to assigning independent and identically distributed $\mathrm{Uniform}(0,1)$ arrival times to all vertices of $G$ and processing them in increasing order of their arrival times (see, e.g., \cite{Coletti2024, Penrose2002, Ritchie2006}). Since blocking occurs only within each connected component, the resulting jamming limit on $G_j$ is determined solely by the relative order of the arrival times within $V(G_j)$. Because these arrival times are independent across distinct components, $N_{G_1}^{\mathrm{seq}}, \dots, N_{G_r}^{\mathrm{seq}}$ are independent, and the convolution formulas follow.
\end{proof}

We now recall a classical result on stochastic orders (Theorem~1.A.3 in \cite{shaked-2007}).

\begin{theorem}[Theorem~1.A.3 of \cite{shaked-2007}]\label{thm:shaked}
Let $X_1, \dots, X_r$ be a set of independent random variables and let $Y_1, \dots, Y_r$ be another set of independent random variables. If $X_i$ is stochastically dominated by $Y_i$ for every $i = 1, \dots, r$, then, for every increasing function $\psi : \mathbb{R}^r \to \mathbb{R}$, $\psi(X_1, \dots, X_r)$ is stochastically dominated by $\psi(Y_1, \dots, Y_r)$.
\end{theorem}

\begin{corollary}\label{cor:conexidad}
Let $G$ be a finite simple disconnected graph with connected components $G_1, G_2, \dots, G_r$. If on every component $G_i$ the sequential model dominates the static model (respectively, the static model dominates the sequential model), then the same holds on $G$.
\end{corollary}

\begin{proof}
By hypothesis, $N_{G_i}^{\mathrm{stat}} \le_{\mathrm{st}} N_{G_i}^{\mathrm{seq}}$ for each $i$, and the collections $\{N_{G_i}^{\mathrm{stat}}\}_{i=1}^{r}$ and $\{N_{G_i}^{\mathrm{seq}}\}_{i=1}^{r}$ each consist of independent random variables by Proposition~\ref{prop:convolucion}. Applying Theorem~\ref{thm:shaked} with the increasing function $\psi(x_1, \dots, x_r) = \sum_{j=1}^{r} x_j$, we conclude that $\sum_{j=1}^{r} N_{G_j}^{\mathrm{stat}} = N_G^{\mathrm{stat}} \le_{\mathrm{st}} \sum_{j=1}^{r} N_{G_j}^{\mathrm{seq}} = N_G^{\mathrm{seq}}$, so the sequential model dominates the static model on $G$. The case in which the static model dominates the sequential one is analogous.
\end{proof}

\subsection{Mixed components: three examples}
\label{sec:mixtos}
When a disconnected graph $G$ is composed of components for which the dominance goes in different directions, the overall behavior is not uniquely determined. We present explicit examples covering all possible outcomes.

\begin{figure}[H]
\centering
\includegraphics[width=0.55\textwidth]{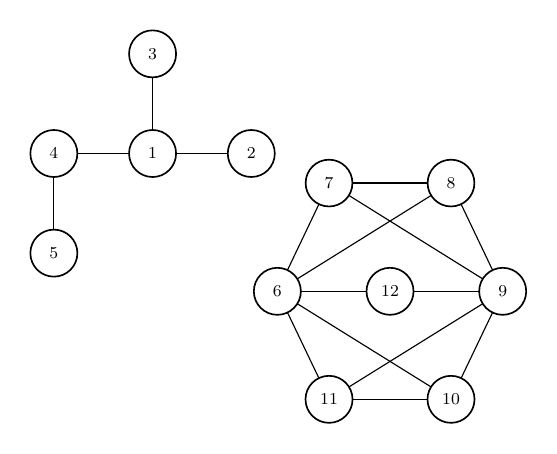}
\caption{Disconnected graph on which the static model dominates the sequential one (Example~1).}
\label{fig:disc_d2}
\end{figure}

\textbf{Example 1 (sequential dominates static $\oplus$ static dominates sequential $\implies$ static dominates sequential).} Consider a disconnected graph $G$ formed by a 5-vertex connected component $G_1$ and a 7-vertex connected component $G_2$ (Figure~\ref{fig:disc_d2}). By Proposition~\ref{thm:pequenos}, the sequential model dominates the static one on the left component $G_1$ (which appears in the fifth row of Table~\ref{tab:subset_grafos_acumulada}, with $\mathbb{P}\left[N_{G_1}^{\mathrm{stat}} = 2\right] = 1/3 > 4/15 = \mathbb{P}\left[N_{G_1}^{\mathrm{seq}} = 2\right]$). The right component $G_2$ has 5 jamming limits, of sizes 2 and 3, with $\mathbb{P}\left[N_{G_2}^{\mathrm{stat}} = 2\right] = 1/5 < 2/7 = \mathbb{P}\left[N_{G_2}^{\mathrm{seq}} = 2\right]$, so on $G_2$ the static model dominates the sequential one. On the combined graph $G$ the static model dominates the sequential one, as shown by the cumulative distributions in Table~\ref{tab:disc_d2}.

\begin{table}[H]
\centering
\renewcommand{\arraystretch}{1.25}
\begin{tabular}{ccc}
\toprule
$k$ & $\mathbb{P}[N_G^{\mathrm{stat}}\leq k]$ & $\mathbb{P}[N_G^{\mathrm{seq}} \leq k]$\\
\midrule
$4$ & $1/15$ & $8/105$\\
$5$ & $7/15$ & $50/105$\\
$6$ & $1$ & $1$\\
\bottomrule
\end{tabular}
\caption{Cumulative distributions of the number of occupied vertices for the disconnected graph of Example~1.}
\label{tab:disc_d2}
\end{table}

\begin{figure}[H]
\centering
\includegraphics[width=0.45\textwidth]{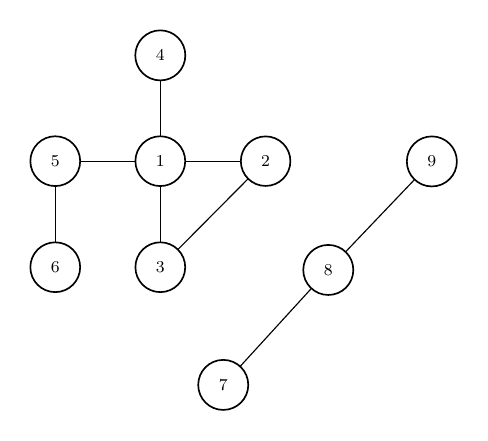}
\caption{Disconnected graph on which the sequential model dominates the static one (Example~2).}
\label{fig:disc_d1}
\end{figure}

\textbf{Example 2 (sequential dominates static $\oplus$ static dominates sequential $\implies$ sequential dominates static).} Consider a disconnected graph $G$ formed by a 3-vertex path $G_1$, on which the sequential model strictly dominates the static one, and the maximum-degree-4 graph $G_2$ of Section~\ref{sec:grado4}, on which the static model strictly dominates the sequential one (Figure~\ref{fig:disc_d1}). On the combined graph $G$ the sequential model dominates the static one, as shown by the cumulative distributions in Table~\ref{tab:disc_d1}.

\begin{table}[H]
\centering
\renewcommand{\arraystretch}{1.25}
\begin{tabular}{ccc}
\toprule
$k$ & $\mathbb{P}[N_G^{\mathrm{stat}}\leq k]$ & $\mathbb{P}[N_G^{\mathrm{seq}}\leq k]$\\
\midrule
$3$ & $1/10$ & $5/72$\\
$4$ & $3/5$ & $17/36$\\
$5$ & $1$ & $1$\\
\bottomrule
\end{tabular}
\caption{Cumulative distributions of the number of occupied vertices for the disconnected graph of Example~2.}
\label{tab:disc_d1}
\end{table}

\begin{figure}[H]
\centering
\includegraphics[width=0.55\textwidth]{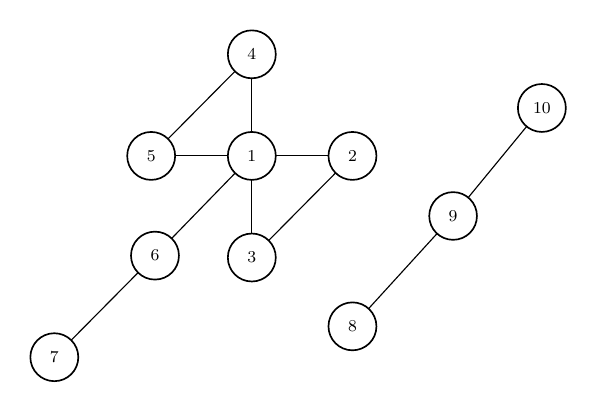}
\caption{Disconnected graph on which neither model dominates the other (Example~3).}
\label{fig:disc_neither}
\end{figure}

\textbf{Example 3 (sequential dominates static $\oplus$ static dominates sequential $\implies$ neither dominates).} Consider a disconnected graph $G$ formed by a 3-vertex path $G_1$ and a 7-vertex component $G_2$ on which the static model dominates the sequential one (Figure~\ref{fig:disc_neither}). The 7-vertex component $G_2$ has 9 jamming limits, of sizes 2 and 3, with $\mathbb{P}\left[N_{G_2}^{\mathrm{stat}} = 2\right] = 1/9 < 6/35 = \mathbb{P}\left[N_{G_2}^{\mathrm{seq}} = 2\right]$, so on $G_2$ the static model strictly dominates the sequential one. On the combined graph $G$ neither model dominates the other: as shown in Table~\ref{tab:disc_neither}, at $k=3$ we have $\mathbb{P}[N_G^{\mathrm{stat}} \le 3] = 1/18 < 2/35 = \mathbb{P}[N_G^{\mathrm{seq}} \le 3]$, whereas at $k=4$ we have $\mathbb{P}[N_G^{\mathrm{stat}} \le 4] = 5/9 > 47/105 = \mathbb{P}[N_G^{\mathrm{seq}} \le 4]$. Since the cumulative distributions cross, no first-order stochastic dominance holds.

\begin{table}[H]
\centering
\renewcommand{\arraystretch}{1.25}
\begin{tabular}{ccc}
\toprule
$k$ & $\mathbb{P}[N_G^{\mathrm{stat}}\leq k]$ & $\mathbb{P}[N_G^{\mathrm{seq}}\leq k]$\\
\midrule
$3$ & $1/18$ & $2/35$\\
$4$ & $5/9$ & $47/105$\\
$5$ & $1$ & $1$\\
\bottomrule
\end{tabular}
\caption{Cumulative distributions of the number of occupied vertices for the disconnected graph of Example~3.}
\label{tab:disc_neither}
\end{table}

These examples illustrate the structural complexity of disconnected graphs when components exhibit opposing stochastic dominance relations.

\section*{Statements and Declarations}

\noindent\textbf{Funding.} Robinson Higuita-Diaz was supported by Universidad Pedag\'ogica y Tecnol\'ogica de Colombia, grant no.~SGI~4141. Alejandro Rold\'an-Correa was supported by Universidad de Antioquia, project no.~2026-89150.

\medskip
\noindent\textbf{Competing interests.} The authors have no competing interests to declare that are relevant to the content of this article.

\medskip
\noindent\textbf{Data and code availability.} No datasets were generated or analysed in this study. The code used for all exact computations reported in this article is publicly available in the GitHub repository \cite{LopezRivera2026Repositorio}.

\medskip
\noindent\rev{\textbf{Acknowledgments.} The authors are grateful to David Jiang and Akash Jangid for their interest in this work and for communicating to us the counterexamples and results presented in Section~\ref{sec:grado3}.}


\begin{thebibliography}{99}

\bibitem{Baram1992} A. Baram and D. Kutasov, Random sequential adsorption on a quasi-one-dimensional lattice: an exact solution, J. Phys. A, {\bf 25}, no. 8, (1992), L493--L498.

\bibitem{Baram1994} A. Baram and D. Kutasov, Random sequential adsorption on a $3 \times \infty$ lattice: an exact solution, J. Phys. A, {\bf 27}, no. 11, (1994), 3683--3687.

\bibitem{Baxter1982} R. J. Baxter, Exactly Solved Models in Statistical Mechanics, Academic Press, London, (1982).

\bibitem{Brightwell2017} G. Brightwell, S. Janson and M. Luczak, The greedy independent set in a random graph with given degrees, Random Structures \& Algorithms, {\bf 51}, no. 4, (2017), 565--586.

\bibitem{Dinning-table} H.-H. Chern, H.-K. Hwang and T.-H. Tsai, Random unfriendly seating arrangement in a dining table, Advances in Applied Mathematics, {\bf 65}, (2015), 38--64.

\bibitem{Coletti2024} C. F. Coletti, S. Gallo, A. Rold\'an-Correa and L. A. Valencia, Fluctuations of the occupation density for a parking process, Journal of Statistical Physics, {\bf 191}, no. 11, (2024), 146.

\bibitem{Dehling2008} H. G. Dehling, S. R. Fleurke and C. K\"ulske, Parking on a random tree, Journal of Statistical Physics, {\bf 133}, no. 1, (2008), 151--157.

\bibitem{Dhara2016} S. Dhara, J. S. H. van Leeuwaarden and D. Mukherjee, Generalized random sequential adsorption on Erd\H{o}s--R\'enyi random graphs, Journal of Statistical Physics, {\bf 164}, (2016), 1217--1232.

\bibitem{Diestel2017} R. Diestel, Graph Theory, 5th ed., Graduate Texts in Mathematics, vol. 173, Springer, Berlin, (2017).

\bibitem{Doslic2024} T. Do\v{s}li\'c, M. Puljiz, S. \v{S}ebek and J. \v{Z}ubrini\'c, Complexity function of jammed configurations of Rydberg atoms, Ars Mathematica Contemporanea, {\bf 24}, no. 4, (2024), \#P4.05.

\bibitem{Evans1993} J. W. Evans, Random and cooperative sequential adsorption, Rev. Mod. Phys., {\bf 65}, no. 4, (1993), 1281--1329.

\bibitem{Fan1992} Y. Fan and J. K. Percus, Random sequential adsorption on a ladder, Journal of Statistical Physics, {\bf 66}, no. 1--2, (1992), 263--271.

\bibitem{Flory1939} P. J. Flory, Intramolecular reaction between neighboring substituents of vinyl polymers, J. Am. Chem. Soc., {\bf 61}, no. 6, (1939), 1518--1521.

\bibitem{Freedman1962} D. Freedman and L. Shepp, An unfriendly seating arrangement, SIAM Rev., {\bf 4}, no. 2, (1962), 150.

\bibitem{Friedman1964} H. D. Friedman, D. Rothman and J. K. Mackenzie, Solution to: An unfriendly seating arrangement (problem 62-3), SIAM Rev., {\bf 6}, no. 2, (1964), 180--182.

\bibitem{theather} K. Georgiou, E. Kranakis and D. Krizanc, Random maximal independent sets and the unfriendly theater seating arrangement problem, Discrete Mathematics, {\bf 309}, no. 16, (2009), 5120--5129.

\bibitem{Higuita2025} R. Higuita-Diaz, L. L\'opez and A. Rold\'an-Correa, Combinatorial aspects of a discrete parking process, Statistics \& Probability Letters, {\bf 233}, (2026), 110682.


\bibitem{Higuita2026Corrigendum} R. Higuita-Diaz, L. L\'opez and A. Rold\'an-Correa, Corrigendum to ``Combinatorial aspects of a discrete parking process'' [Statist. Probab. Lett. 233 (2026) 110682], Statistics \& Probability Letters, (2026), 111009. \url{https://doi.org/10.1016/j.spl.2026.111009}


\bibitem{Krapivsky2020} P. L. Krapivsky, Large deviations in one-dimensional random sequential adsorption, Physical Review E, {\bf 102}, (2020), 062108.

\bibitem{KrapivskyLuck2023} P. L. Krapivsky and J. M. Luck, A renewal approach to configurational entropy in one dimension, J. Phys. A: Math. Theor., {\bf 56}, no. 25, (2023), 255001.

\bibitem{Krivelevich2024} M. Krivelevich, T. M\'esz\'aros, P. Michaeli and C. Shikhelman, Greedy maximal independent sets via local limits, Random Structures \& Algorithms, {\bf 64}, no. 4, (2024), 986--1015.

\bibitem{LopezRivera2026Repositorio} L. L\'opez, Stochastic Dominance DPP Graphs, GitHub repository, (2026), \url{https://github.com/lykoslopezrivera-blip/Stochastic_Dominance_DPP_Graphs}.

\bibitem{OEIS_A001349} OEIS Foundation Inc., Number of connected graphs with n unlabeled nodes, The On-Line Encyclopedia of Integer Sequences, sequence A001349, (2026), \url{https://oeis.org/A001349}.

\bibitem{Oh2017} S. Oh, Maximal independent sets on a grid graph, Discrete Mathematics, {\bf 340}, no. 12, (2017), 2762--2768.

\bibitem{Page1959} E. S. Page, The distribution of vacancies on a line, J. R. Stat. Soc. Ser. B, {\bf 21}, no. 2, (1959), 364--374.

\bibitem{Penrose2002} M. D. Penrose, Limit theorems for monotonic particle systems and sequential deposition, Stochastic Process. Appl., {\bf 98}, no. 2, (2002), 175--197.

\bibitem{PenroseSudbury2005} M. D. Penrose and A. Sudbury, Exact and approximate results for deposition and annihilation processes on graphs, Ann. Appl. Probab., {\bf 15}, no. 1B, (2005), 853--889.

\bibitem{Renyi1958} A. R\'enyi, On a one-dimensional problem concerning random space filling, Magy. Tud. Akad. Mat. Kutat\'o Int. K\"ozl., {\bf 3}, no. 1--2, (1958), 109--127.

\bibitem{Ritchie2006} T. L. Ritchie, Construction of the thermodynamic jamming limit for the parking process and other exclusion schemes on $\mathbb{Z}^d$, Journal of Statistical Physics, {\bf 122}, no. 3, (2006), 381--398.

\bibitem{shaked-2007} M. Shaked and J. G. Shanthikumar, Stochastic Orders, Springer, New York, (2007).

\bibitem{Shneer2017} S. Shneer and P. M. van de Ven, Per-site occupancy in the discrete parking problem, Statistics \& Probability Letters, {\bf 120}, (2017), 141--146.

\bibitem{SpahnZeilberger2024} G. Spahn and D. Zeilberger, Enumerating seating arrangements that obey social distancing, Journal of Symbolic Computation, {\bf 126}, (2025), 102344.

\bibitem{Sudbury2009} A. Sudbury, Random sequential adsorption on random trees, Journal of Statistical Physics, {\bf 136}, no. 1, (2009), 51--58.

\bibitem{Talbot2000} J. Talbot, G. Tarjus, P. R. Van Tassel and P. Viot, From car parking to protein adsorption: an overview of sequential adsorption processes, Colloids Surf. A, {\bf 165}, no. 1, (2000), 287--324.

\end{thebibliography}
\end{document}